\documentclass[12pt,a4paper]{article}

\usepackage{amssymb}

\usepackage{amsmath}

\newtheorem{theorem}{Theorem}[section]
\newtheorem{definition}[theorem]{Definition}
\newtheorem{lemma}[theorem]{Lemma}
\newtheorem{corollary}[theorem]{Corollary}
\newtheorem{proposition}[theorem]{Proposition}
\newtheorem{example}[theorem]{Example}

\begin{document}
\title{Gr\"{o}bner-Shirshov bases for dialgebras\footnote{Supported by the
NNSF of China (Nos.10771077, 10911120389) and the NSF of Guangdong
Province (No.06025062).}}
\author{
L. A. Bokut\footnote {Supported by RFBR 01-09-00157, LSS--344.2008.1
and SB RAS Integration grant No. 2009.97 (Russia).} \\
{\small \ School of Mathematical Sciences, South China Normal
University}\\
{\small Guangzhou 510631, P. R. China}\\
{\small Sobolev Institute of Mathematics, Russian Academy of
Sciences}\\
{\small Siberian Branch, Novosibirsk 630090, Russia}\\
{\small Email: bokut@math.nsc.ru}\\
\\
 Yuqun
Chen\footnote {Corresponding author.} \  and Cihua Liu\\
{\small \ School of Mathematical Sciences, South China Normal
University}\\
{\small Guangzhou 510631, P. R. China}\\
{\small Email: yqchen@scnu.edu.cn}\\
{\small langhua01duo@yahoo.com.cn}}

\date{}

\maketitle \noindent\textbf{Abstract:} In this paper, we define the
Gr\"{o}bner-Shirshov basis for a dialgebra. The Composition-Diamond
lemma for dialgebras is given then. As results, we give
Gr\"{o}bner-Shirshov bases for the universal enveloping algebra of a
Leibniz algebra, the bar extension of a dialgebra, the free product
of two dialgebras, and Clifford dialgebra. We obtain some normal
forms for  algebras mentioned the above.



\section{Introduction}

J.-L. Loday  (1995, \cite{Lo95}) gave the definition of a new class
of algebras, dialgebras, which is closely connected to his notion of
Leibniz algebras (1993, \cite{Lo93}) in the same way as associative
algebras connected to Lie algebras. In the manuscript \cite{Lo99},
J.-L. Loday found a normal form of elements of a free dialgebra.
Here we continue to study free dialgebras and prove the
Composition-Diamond lemma for dialgebras. As it is well known, this
kind of lemma is the cornerstone of the theory of Gr\"{o}bner and
Gr\"{o}bner-Shirshov bases (see, for example, \cite{bo05} and cited
literature). In commutative-associative case, this lemma is
equivalent to the Main Buchberger's Theorem (\cite{bu65, bu70}). For
Lie and associative algebras, this is the Shirshov's lemma
\cite{sh62} (see also L.A. Bokut \cite{bo72, bo76}, G. Bergman
\cite{be78},  L.A. Bokut and Y. Chen \cite{bc1}). As results, we
obtain  Gr\"{o}bner-Shirshov bases for the universal enveloping
algebra of a Leibniz algebra, the bar extension of a dialgebra, the
free product of two dialgebras, and Clifford dialgebra. By using our
Composition-Diamond lemma for dialgebras (Theorem \ref{cd}), we
obtain some normal forms for  algebras mentioned the above.
Moreover, we get another proof of the M. Aymon, P.-P. Grivel's
result (\cite{A}) on the Poincare-Birkhoff-Witt theorem for Leibniz
algebras  (see P. Kolesnikov \cite{K} for other proof).

\section{Preliminaries}

\begin{definition}\label{l2}
Let $k$ be a field. A $k$-linear space $D$ equipped with two
bilinear multiplications $\vdash$  and $\dashv $  is called a
dialgebra, if both $\vdash$ and $\dashv$ are associative and
\begin{eqnarray*}
a\dashv(b\vdash c)&=&a\dashv b\dashv c \\
(a\dashv b)\vdash c&=&a\vdash b\vdash c \\
a\vdash(b\dashv c)&=&(a\vdash b)\dashv c
\end{eqnarray*}
for any $a, \ b, \ c\in D$.
\end{definition}

\begin{definition}
Let $D$ be a dialgebra, $B\subset D$. Let us define diwords
 of $D$ in the set $B$ by induction:
\begin{enumerate}
\item[(i)]\  $b=(b), \ b\in B$ is a diword in $B$ of length $|b|=1$.
\item[(ii)]\  $(u)$ is called a diword in $B$ of length $|(u)|=n$, if
$(u)=((v)\dashv(w))$ or $(u)=((v)\vdash(w))$, where $(v), \ (w)$ are
diwords in $B$ of length $k, \ l$ respectively and $k+l=n$.
\end{enumerate}
\end{definition}

\begin{proposition}\label{1}(\cite{Lo99})
Let $D$ be a dialgebra and $B\subset D$. Any diword of $D$ in the
set $B$ is equal to a diword in $B$ of the form

\begin{equation}\label{(1)}
(u)=b_{-m}\vdash \cdots \vdash b_{-1} \vdash b_0 \dashv b_1 \dashv
\cdots \dashv b_n
\end{equation}
where $b_i \in B, \ -m\leq i\leq n, \ m\geq 0, \ n\geq 0$. Any
bracketing of the right side of (\ref{(1)}) gives the same result. \
\ $\square$
\end{proposition}

\begin{definition}
Let $X$ be a set. A free dialgebra $D(X)$ generated by $X$ over $k$
is defined in a usual way by the following commutative diagram:

\setlength {\unitlength}{1cm}
\begin{picture}(7, 3)
\put(4.2,2.3){\vector(1,0){1.7}} \put(4.1, 2.0){\vector(0,-1){1.3}}
\put(5.9,2.1){\vector(-1,-1){1.6}} \put(3.9,0.2){$D$}
\put(6,2.2){$D(X)$} \put(3.9,2.2){$X$} \put(4.9, 2.4){$i$}
\put(3.4,1.3){$\forall\varphi$} \put(5.3,1){$\exists!\varphi^*$
\mbox{ (homomorphism)}}
\end{picture}\\
where $D$ is any dialgebra.
\end{definition}

In \cite{Lo99}, a construction of a free dialgebra is given.

\begin{proposition}(\cite{Lo99})
Let $D(X)$ be a free dialgebra over $k$ generated by $X$. Any diword
in $D(X)$ is equal to the unique diword of the form
\begin{equation}\label{(2)}
[u]=x_{-m}\vdash \cdots \vdash x_{-1} \vdash x_0 \dashv x_1 \dashv
\cdots \dashv x_n \triangleq x_{-m} \cdots  x_{-1} \dot{x_0} x_1
\cdots x_n
\end{equation}
where $x_i \in X, \ m\geq 0, \ n\geq 0$, and $x_0$ is called the
center of the normal diword $[u]$. We call $[u]$ a normal diword (in
$X$) with the associative word $u, u\in X^*$. Clearly, if $[u]=[v]$,
then $u=v$. In (\ref{(2)}). Let $[u], \ [v]$ be two normal diwords.
Then $[u]\vdash[v]$ is the normal diword $[uv]$ with the center at
the center of $[v]$. Accordingly, $[u]\dashv[v]$ is the normal
diword $[uv]$ with the center at the center of $[u]$. \ \ $\square$
\end{proposition}

\begin{example}
$$
(x_{-1}\vdash x_0\dashv x_1)\vdash(y_{-1}\vdash y_0\dashv
y_1)=x_{-1}\vdash x_0\vdash x_1\vdash y_{-1}\vdash y_0\dashv y_1,
$$
$$
(x_{-1}\vdash x_0\dashv x_1)\dashv(y_{-1}\vdash y_0\dashv y_1)=
x_{-1}\vdash x_0\dashv x_1\dashv y_{-1}\dashv y_0\dashv y_1. \ \
\square
$$
\end{example}

\section{Composition-Diamond lemma for dialgebras }

Let $X$ be a well ordered set, $D(X)$ the free dialgebra over $k$,
$X^*$ the free monoid generated by $X$ and  $[X^*]$ the set of
normal diwords in $X$. Let us define the deg-lex ordering on $[X^*]$
in the following way: for any $[u],[v]\in [X^*]$,
$$
[u]<[v] \Longleftrightarrow wt([u])<wt([v]) \ \
\mbox{lexicographicaly},
$$
where
$$
wt([u])=(n+m+1,m, x_{-m}, \cdots,x_0, \cdots, x_n)
$$
if $[u]=x_{-m} \cdots  x_{-1} \dot{x_0} x_1  \cdots x_n$.

Throughout the paper, we will use this ordering.

It is easy to see that the ordering $<$ is satisfied the following
properties:
$$
[u]<[v]\Longrightarrow x\vdash[u]<x\vdash[v], \ [u]\dashv x<
[v]\dashv x, \ \mbox{for any } \ x\in X.
$$

Any polynomial $f\in D(X)$ has the form
$$
f=\sum_{[u]\in
[X^*]}f([u])[u]=\alpha[\overline{f}]+\sum{\alpha}_i[u_i],
$$
where $[\overline{f}], \ [u_i]$ are normal diwords in $X$,
$[\overline{f}]>[u_i], \ \alpha , \ {\alpha}_i, \ f([u])\in k,\
\alpha\neq 0$. We call $[\overline{f}]$ the leading term of $f$.
Denote $suppf$ by the set $\{ [u]|f([u])\neq 0\}$ and $deg(f)$ by
$|[\overline{f}]|$. $f$ is called monic if $\alpha=1$. $f$ is called
left (right) normed if $f=\sum{\alpha}_i u_i \dot{x_i} \ \
(f=\sum{\alpha}_i \dot{x_i} u_i)$, where each ${\alpha}_i\in k, \
x_i\in X$ and $u_i\in X^*$.

If $[u], \ [v]$ are both left normed or both right normed, then it
is clear that for any $[w]\in [X^*]$,
\begin{eqnarray*}
[u]<[v]\Longrightarrow &&[u]\vdash [w]<[v]\vdash [w], \
[w]\vdash[u]<[w]\vdash[v],\\
&& [u]\dashv [w]<[v]\dashv [w], \ [w]\dashv[u]<[w]\dashv[v].
\end{eqnarray*}

Let $S\subset D(X)$. By an $S$-diword $g$ we will mean  a diword in
$\{X\cup S\}$ with only one occurrence of $s\in S$. If this is the
case and $g=(asb)$ for some $a,b\in X^*, \ s\in S$, we also call $g$
an $s$-diword.

From Proposition \ref{1} it follows that any $s$-diword is equal to
\begin{equation}\label{(3)}
[asb]=x_{-m}\vdash \cdots \vdash x_{-1} \vdash x_0 \dashv x_1 \dashv
\cdots \dashv x_n|_{x_k\mapsto s}
\end{equation}
where $-m\leq k\leq n,  \ s\in S, \ x_i\in X, \ -m\leq i\leq n$. To
be more precise, $[asb]=[a\dot{s}b]$ if $k=0$;
$[asb]=[asb_1\dot{x_0}b_2]$ if $k<0$ and $[asb]=[a_1\dot{x_0}a_2sb]$
if $k>0$. If the center of the $s$-diword $[asb]$ is in $a$, then we
denote it by $[\dot{a}sb]=[a_1\dot{x_0}a_2sb]$. Similarly,
$[as\dot{b}]=[asb_1\dot{x_0}b_2]$ (of course, either $a_i$ or $b_i$
may be empty).

\begin{definition}
The $s$-diword (\ref{(3)}) is called a normal $s$-diword if one of
the following conditions holds:
\begin{enumerate}
\item[(i)] \ $k=0$,
\item[(ii)] \ $k<0$ and $s$ is left normed,
\item[(iii)] \ $k>0$ and $s$ is right normed.
\end{enumerate}
We call a normal $s$-diword $[asb]$ a left (right) normed $s$-diword
if both $s$ and $[asb]$ are left (right) normed. In particulary, $s$
is a left (right) normed $s$-diword if $s$ is left (right) normed
polynomial.
\end{definition}

The following lemma follows from the above properties of the
ordering $<$.
\begin{lemma}
For a normal $s$-diword $[asb]$, the leading term of $[asb]$ is
equal to $[a[\overline{s}]b]$, that is,
$\overline{[asb]}=[a[\overline{s}]b]$. More specifically, if
$$
[asb]=x_{-m}\vdash \cdots \vdash x_{-1} \vdash x_0 \dashv x_1 \dashv
\cdots \dashv x_n|_{x_k\mapsto s},
$$
then corresponding to $k=0,\ k<0,\ k>0$, respectively, we have
$$
\overline{x_{-m}\vdash \cdots \vdash x_{-1} \vdash s \dashv x_1
\dashv \cdots \dashv x_n}=x_{-m}\vdash \cdots \vdash x_{-1} \vdash
[\overline{s}] \dashv x_1 \dashv \cdots \dashv x_n,
$$
$$
\overline{x_{-m}\vdash \cdots \vdash s\vdash \cdots \vdash x_0
 \dashv \cdots \dashv x_n}=x_{-m}\vdash \cdots \vdash
[\overline{s}]\vdash \cdots \vdash x_0 \dashv \cdots \dashv x_n,
$$
$$
\overline{x_{-m}\vdash \cdots \vdash x_0 \dashv \cdots \dashv s
\dashv \cdots \dashv x_n}={x_{-m}\vdash \cdots \vdash x_0 \dashv
\cdots \dashv [\overline{s}] \dashv \cdots \dashv x_n}. \ \ \square
$$
\end{lemma}

\ \

Now, we define compositions of polynomials in $D(X)$.

\begin{definition}\label{2.9}
Let the ordering $<$ be as before and $f,g\in D(X)$ with $f,g$
monic.
\begin{enumerate}
\item[1)] \ Composition of left (right) multiplication.

Let $f$ be not a right normed polynomial and $x\in X$. Then $x\dashv
f$ is called the composition of left multiplication. Clearly,
$x\dashv f$ is a right normed polynomial (or 0).

Let $f$ be not a left normed polynomial and $x\in X$. Then $f\vdash
x$ is called the composition of right multiplication. Clearly,
$f\vdash x$ is  a left normed polynomial (or 0).

\item[2)] \ Composition of inclusion.

Let
$$
[w]=[\overline{f}]=[a [\overline{g}] b],
$$
where $[agb]$ is a normal $g$-diword. Then
$$
(f,g)_{[w]}=f-[agb]
$$
is called the composition of inclusion. The transformation $f\mapsto
f-[agb]$ is called the elimination of leading diword (ELW) of $g$ in
$f$, and $[w]$ is called the ambiguity of $f$ and $g$.

\item[3)] \ Composition of intersection.

Let
$$[w]=[[\overline{f}]b]=[a[\overline{g}]], \ |\overline{f}|+|\overline{g}|>|w|,$$
where $[fb]$ is a normal $f$-diword and $[ag]$ a normal $g$-diword.
Then
$$(f,g)_{[w]}=[fb]-[ag]$$
is called the composition of intersection, and $[w]$ is called the
ambiguity of $f$ and $g$.
\end{enumerate}
\end{definition}

\noindent \textbf{Remark} \ In the Definition \ref{2.9}, for the
case of 2) or 3), we have $ \overline{(f,g)_{[w]}}<[w]. $ For the
case of 1),  $deg(x\dashv f)\leq deg(f)+1$ and $deg(f\vdash x)\leq
deg(f)+1$.

\begin{definition}
Let the ordering $<$ be as before, $S\subset D(X)$ a monic set and
$f,g\in S$.
\begin{enumerate}
\item[1)] \ Let $x\dashv f$ be a composition of left
multiplication. Then $x\dashv f$ is called trivial  modulo $S$,
denoted by \ \ $x\dashv f\equiv0 \ mod(S)$, \ \ if
$$
x\dashv f=\sum{\alpha}_i[a_i s_i b_i],
$$
where each $\alpha_i\in k, \ a_i,b_i\in X^*, \ s_i\in S, \ [a_i s_i
b_i]$ right normed $s_i$-diword and $|[a_i [\overline{s_i}]
b_i]|\leq deg(x\dashv f)$.

\ \

Let $f\vdash x$ be a composition of right multiplication. Then
$f\vdash x$ is called trivial  modulo $S$, denoted by \ \ $f\vdash
x\equiv0 \ mod(S)$,  \ \ if
$$
f\vdash x=\sum{\alpha}_i[a_i s_i b_i],
$$
where each $\alpha_i\in k, \ a_i,b_i\in X^*, \ s_i\in S, \ [a_i s_i
b_i]$ left normed  $s_i$-diword and $|[a_i [\overline{s_i}]
b_i]|\leq deg(f\vdash x)$.

\item[2)]  \ Composition $(f,g)_{[w]}$ of inclusion (intersection) is
called trivial modulo $(S,[w])$, denoted by  \ \ $(f,g)_{[w]}\equiv0
\ mod(S,[w])$, \ \  if
$$
(f,g)_{[w]}=\sum{\alpha}_i[a_i s_i b_i],
$$
where each $\alpha_i\in k, \ a_i,b_i\in X^*, \ s_i\in S, \ [a_i s_i
b_i]$ normal $s_i$-diword, $[a_i [\overline{s_i}] b_i]<[w]$ and each
$[a_i s_i b_i]$ is right (left) normed  $s_i$-diword whenever either
both $f$ and $[agb]$ or both $[fb]$ and $[ag]$ are right (left)
normed $S$-diwords.
\end{enumerate}

We call the set $S$ a Gr\"{o}bner-Shirshov basis in $D(X)$ if any
composition of polynomials in $S$ is trivial modulo $S$ (and $[w])$.
\end{definition}

\ \

The following lemmas play key role in the proof of Theorem \ref{cd}.

\begin{lemma}\label{l4}
Let $S\subset D(X)$ and $[asb]$ an
 $s$-diword, $s\in S$. Assume that each composition of right and left
multiplication is trivial modulo $S$. Then, $[asb]$ has a
presentation:
$$[asb]=\sum\alpha_i[a_is_ib_i],
$$
where each $\alpha_i\in k, \ s_i\in S, \ a_i,b_i\in X^*$ and each
$[a_i s_i b_i]$ is normal $s_i$-diword.
\end{lemma}
\textbf{Proof.} \ Following Proposition \ref{1}, we assume  that
$$
[asb]=x_{-m}\vdash \cdots \vdash x_{-1} \vdash x_0 \dashv x_1 \dashv
\cdots \dashv x_n|_{x_k\mapsto s}.
$$
 There are three cases to consider.

Case 1. $k=0$. Then $[asb]$ is a normal $s$-diword.

Case 2. $k<0$. Then $[asb]=a\vdash (s\vdash x_{k+1})\vdash b, k<-1$
or $[asb]=a\vdash (s\vdash x_{0})\dashv b$. If $s$ is left normed
then $[asb]$ is a normal $s$-diword.  If $s$ is not left normed then
for the composition $s\vdash x_{k+1} \ \ (k<0)$ of right
multiplication, we have
$$
s\vdash x_{k+1}=\sum{\alpha}_i[a_i s_i b_i],
$$
where each $\alpha_i\in k, \ a_i,b_i\in X^*, \ s_i\in S$ and $[a_i
s_i b_i]$ is left normed  $s_i$-diword. Then
$$
[asb]=\sum{\alpha}_i(a\vdash[a_i s_i b_i]\vdash b)
$$
or
$$
[asb]=\sum{\alpha}_i(a\vdash[a_i s_i b_i]\dashv b)
$$
 is a  linear combination of normal $s_i$-diwords.

Case 3. $k>0$ is similar to the Case 2.  \ \ $\square$

\begin{lemma}\label{l5}
Let $S\subset D(X)$ and  each composition $(f,g)_{[w]}$ in $S$ of
inclusion (intersection)  trivial modulo $(S,[w])$. Let
$[a_1s_1b_1]$ and $[a_2s_2b_2]$ be normal $S$-diwords such that
$[w]=[a_1[\bar{s_1}]b_1]=[a_2[\bar{s_2}]b_2]$, where $s_1,s_2\in S,\
a_1,a_2,b_1,b_2\in X^*$. Then,
$$
[a_1s_1b_1]\equiv[a_2s_2b_2] \ \ mod(S,[w]),
$$
i.e., $[a_1s_1b_1]-[a_2s_2b_2]=\sum{\alpha}_i[a_i s_i b_i], $ where
each $\alpha_i\in k, \ a_i,b_i\in X^*, \ s_i\in S, \ [a_i s_i b_i]$
normal $s_i$-diword and $[a_i [\overline{s_i}] b_i]<[w]$.
\end{lemma}
\noindent{\bf Proof.} \ In the following, all letters $a,b,c$ with
indexis are words and $s_1,s_2,s_j\in S$.

Because $a_1\bar{s_1}b_1=a_2\bar{s_2}b_2$ as ordinary words, there
are three cases to consider.

Case 1. Subwords $\overline{s_1}, \overline{s_2}$ have empty
intersection. Assume, for example,  that $b_1=b\overline{s_2}b_2$
and $a_2=a_1\overline{s_1}b$. Because any normal $S$-diword may be
bracketing in any way,  we have
$$
[a_2s_2b_2]-[a_1s_1b_1]=(a_1s_1(b(s_2-[\overline{s_2}])b_2))-
((a_1(s_1-[\overline{s_1}])b)s_2b_2).
$$
For any $[t]\in supp(s_2-[\overline{s_2}])$, we prove that
$(a_1s_1b[t]b_2)$ is a normal $s_1$-diword. There are five cases to
consider.

1.1 $[w]=[\dot{a_1}[\overline{s_1}]b[\overline{s_2}]b_2]$;

1.2 $[w]=[a_1\dot{[\overline{s_1}]}b[\overline{s_2}]b_2]$;

1.3 $[w]=[a_1[\overline{s_1}]\dot{b}[\overline{s_2}]b_2]$;

1.4 $[w]=[a_1[\overline{s_1}]b\dot{[\overline{s_2}]}b_2]$;

1.5 $[w]=[a_1[\overline{s_1}]b[\overline{s_2}]\dot{b_2}]$.

For 1.1, since $[a_1s_1b_1]$ and $[a_2s_2b_2]$ are normal
$S$-diwords, both $s_1$ and  $s_2$ are right normed by the
definition, in particular, $[t]$ is right normed. It follows that
$(a_1s_1b[t]b_2)=[\dot{a_1}s_1b[t]b_2]$ is a normal $s_1$-diword.

For 1.2, it is clear that $(a_1s_1b[t]b_2)$ is a normal $s_1$-diword
and $[t]$ is right normed.

For 1.3, 1.4 and 1.5, since $[a_1s_1b_1]$ is normal $s_1$-diword,
$s_1$ is left normed by the definition, which implies that
$(a_1s_1b[t]b_2)$ is a normal $s_1$-diword. Moreover,  $[t]$ is
right normed, if 1.3, and left normed, if 1.5.

Clearly, for all cases, we have
$\overline{[a_1s_1b[t]b_2]}=[a_1[\overline{s_1}]b[t]b_2]<
[a_1[\overline{s_1}]b[\overline{s_2}]b_2]=[w]$.

Similarly, for any $[t]\in supp(s_1-[\overline{s_1}])$,
$(a_1[t]bs_2b_2)$ is a normal $s_2$-diword and
$[a_1[t]b[\overline{s_2}]b_2]<[w]$.

Case 2. Subwords $\overline{s_1}$ and $\overline{s_2}$ have
non-empty intersection $c$. Assume, for example, that $b_1=bb_2, \
a_2=a_1a, \ w_1=\overline{s_1}b=a\overline{s_2}=acb$.

There are following five cases to consider:

2.1 $[w]=[\dot{a_1}[\overline{s_1}]bb_2]$;

2.2 $[w]=[a_1[\overline{s_1}]b\dot{b_2}]$;

2.3 $[w]=[a_1\dot{a}cbb_2]$;

2.4 $[w]=[a_1a\dot{c}bb_2]$;

2.5 $[w]=[a_1ac\dot{b}b_2]$.

 Then
$$
[a_2s_2b_2]-[a_1s_1b_1]=(a_1([as_2]-[s_1b])b_2)=(a_1(s_1,s_2)_{[w_1]}b_2),
$$

where $[w_1]=[acb]=[[\overline{s_1}]b]=[a[\overline{s_2}]]$ is as
follows:

2.1 $[w_1]$ is right normed;

2.2 $[w_1]$ is left normed;

2.3 $[w_1]=[\dot{a}cb]$;

2.4 $[w_1]=[a\dot{c}b]$;

2.5 $[w_1]=[ac\dot{b}]$.

Since each composition $(f,g)_{[w]}$ in $S$ is trivial modulo
$(S,[w])$, there exist $\beta_j\in k, \ u_j,v_j\in X^*, \ s_j\in S$
such that $[s_1b]-[as_2]=\sum_j\beta_j[u_js_jv_j]$, where each
$[u_js_jv_j]$ is normal $S$-diword and
$[u_j[\overline{s_j}]v_j]<[w_1]=[acb]$. Therefore,
$$
[a_2s_2b_2]-[a_1s_1b_1]=\sum_j\beta_j(a_1[u_js_jv_j]b_2).
$$
Now, we prove that each $(a_1[u_js_jv_j]b_2)$ is normal $s_j$-diword
and
$\overline{(a_1[u_js_jv_j]b_2)}<[w]=[a_1[[\overline{s_1}]b]b_2]$.

For 2.1, since $[\dot{a_1}s_1bb_2]$ and $[\dot{a_1}as_2b_2]$  are
normal $S$-diwords, both $[s_1b]$ and $[as_2]$ are right normed
$S$-diwords. Then, by  definition, each $[u_js_jv_j]$ is right
normed $S$-diword, and so each
$(a_1[u_js_jv_j]b_2)=[\dot{a_1}u_js_jv_jb_2]$ is  normal $S$-diword.

For 2.2, both $[s_1b]$ and $[as_2]$
 must be left normed $S$-diwords. Then, by definition, each $[u_js_jv_j]$
is left normed $S$-diword, and so each
$(a_1[u_js_jv_j]b_2)=[a_1u_js_jv_j\dot{b_2}]$ is  normal $S$-diword.

For 2.3, 2.4 or 2.5, by noting that
$(a_1[u_js_jv_j]b_2)=((a_1)\vdash[u_js_jv_j]\dashv(b_2))$ and
$[u_js_jv_j]$ is normal $S$-diword, $(a_1[u_js_jv_j]b_2)$ is also
normal $S$-diword.

Now, for all cases, we have $\overline{[a_1u_js_jv_jb_2]}=
[a_1u_j[\overline{s_j}]v_jb_2]<[w]=[a_1[acb]b_2]$.

Case 3. One of the subwords $\overline{s_1}$ and $\overline{s_2}$
contains another as a subword. Assume, for example, that $b_2=bb_1,
\ a_2=a_1a, \ w_1=\overline{s_1}=a\overline{s_2}b$.

Again there are following five cases to consider:

2.1 $[w]=[\dot{a_1}a[\overline{s_2}]bb_1]$;

2.2 $[w]=[a_1a[\overline{s_2}]b\dot{b_1}]$;

2.3 $[w]=[a_1\dot{a}[\overline{s_2}]bb_1]$;

2.4 $[w]=[a_1a\dot{[\overline{s_2}]}bb_1]$;

2.5 $[w]=[a_1a[\overline{s_2}]\dot{b}b_1]$.

Then
$$
[a_1s_1b_1]-[a_2s_2b_2]=(a_1(s_1-as_2b)b_1)=(a_1(s_1,s_2)_{[w_1]}b_1).
$$
It is similar to the proof of the Case 2 that we have
$[a_1s_1b_1]\equiv[a_2s_2b_2] \ \ mod(S,[w])$.  \ \ $\square$

\ \
\begin{definition}
Let $S\subset D(X)$. Then
$$
Irr(S)\triangleq \{u\in [X^*]|u\neq [a[\overline{s}]b],s\in S,a,b\in
X^*, \ [asb] \mbox{ is normal s-diword}\}.
$$
\end{definition}

\ \

\begin{lemma}\label{e4}
Let $S\subset D(X)$ and $h\in D(X)$. Then $h$ has a representation
\begin{equation*}
h=\sum_{I_1}\alpha_i[u_i]+\sum_{I_2}\beta_j[a_js_jb_j]
\end{equation*}
where $[u_i]\in Irr(S), \ i\in I_1, [a_js_jb_j]$ normal
$s_j$-diwords, $s_j\in S,\ j\in I_2$ with
$[a_1[\overline{s_1}]b_1]>[a_2[\overline{s_2}]b_2]>\cdots>[a_n[\overline{s_n}]b_n]$.
\end{lemma}
\noindent{\bf Proof.} \ Let $h=\alpha_1[\overline{h}]+\cdots$. We
prove the result by induction on $[\overline{h}]$.

If $[\overline{h}]\in Irr(S)$, then take $[u_1]=[\overline{h}]$ and
$h_1=h-\alpha_1[u_1]$. Clearly, $[\overline{h_1}]<[\overline{h}]$ or
$h_1=0$.

If $[\overline{h}]\not\in Irr(S)$, then
$[\overline{h}]=[a_1[\overline{s_1}]b_1]$ with $[a_1s_1b_1]$ a
normal $s_1$-diword. Let $h_1=h-\beta_1[a_1s_1b_1]$. Then
$[\overline{h_1}]<[\overline{h}]$ or $h_1=0$.  \ \ $\square$

\ \

The following theorem is the main result.

\begin{theorem}\label{cd}
(Composition-Diamond lemma) \ Let $S\subset D(X)$ be a monic set and
the ordering $<$ as before, $Id(S)$ is the ideal generated by $S$.
Then $(i)\Rightarrow (ii)\Leftrightarrow (ii)'\Leftrightarrow
(iii)$, where
\begin{enumerate}
\item[(i)] \ $S$ is a Gr\"{o}bner-Shirshov basis in $D(X)$.
\item[(ii)] \
$f\in Id(S)\Rightarrow [\overline{f}]=[a[\overline{s}]b]$ for some
$s\in S, \ a,b\in X^*$ and $[asb]$ a normal $S$-diword.
\item[$(ii)'$]
$f\in Id(S)\Rightarrow
f=\alpha_1[a_1s_1b_1]+\alpha_2[a_2s_2b_2]+\cdots+\alpha_n[a_ns_nb_n]$
with
$[a_1[\overline{s_1}]b_1]>[a_2[\overline{s_2}]b_2]>\cdots>[a_n[\overline{s_n}]b_n],
$ where $[a_is_ib_i]$ is normal $s_i$-diword, $i=1,2,\cdots,n$.
\item[(iii)] \
The set $ Irr(S)$ is a linear basis of the dialgebra
$D(X|S)=D(X)/Id(S)$ generated by $X$ with defining relations $S$.
\end{enumerate}
\end{theorem}
\noindent{\bf Proof.} \ $(i)\Rightarrow (ii)$. \ Let $S$ be a
Gr\"{o}bner-Shirshov basis and $0\neq f\in Id(S)$. We may assume, by
Lemma \ref{l4}, that
$$
f=\sum_{i=1}^n\alpha_i[a_is_ib_i],
$$
where each $\alpha_i\in k, \ a_i,b_i\in X^*, \ s_i\in S$ and $
[a_is_ib_i]$ normal $S$-diword. Let
$$
[w_i]=[a_i[\overline{s_i}]b_i], \
[w_1]=[w_2]=\cdots=[w_l]>[w_{l+1}]\geq\cdots,\ l\geq 1.
$$
We will use  induction on $l$ and $[w_1]$ to prove that
$[\overline{f}]=[a[\overline{s}]b]$ for some $s\in S \ \mbox{and} \
a,b\in X^*$. If $l=1$, then
$[\overline{f}]=\overline{[a_1s_1b_1]}=[a_1[\overline{s_1}]b_1]$ and
hence the result holds. Assume that $l\geq 2$. Then, by Lemma
\ref{l5}, we have $[a_1s_1b_1]\equiv[a_2s_2b_2] \ \ mod(S,[w_1])$.

Thus, if $\alpha_1+\alpha_2\neq 0$ or $l>2$, then the result follows
from  induction on $l$. For the case $\alpha_1+\alpha_2= 0$ and
$l=2$, we use  induction on $[w_1]$. Now, the result follows.

$(ii)\Rightarrow (ii)'$. \ Assume (ii) and $0\neq f\in Id(S)$. Let
$f=\alpha_1[\overline{f}]+\sum_{[u_i]<[\overline{f}]}\alpha_i[u_i]$.
Then, by (ii), $[\overline{f}]=[a_1[\overline{s_1}]b_1]$, where
$[a_1s_1b_1]$ is a normal $S$-diword. Therefore,
$$
f_1=f-\alpha_1[a_1s_1b_1], \ [\overline{f_1}]<[\overline{f}] \mbox{
or } \ f_1=0, \ f_1\in Id(S).
$$
 Now, by using
induction on $[\overline{f}]$, we have $(ii)'$.

$(ii)'\Rightarrow (ii)$. This part is clear.

 $(ii)\Rightarrow
(iii)$. \  Assume $(ii)$. Then by Lemma \ref{e4}, $Irr(S)$ spans
$D(X|S)$ as $k$-space.

Suppose that $0\neq\sum\alpha_i[u_i]\in Id(S)$ where
$[u_1]>[u_2]>\cdots, \ [u_i]\in Irr(S)$.  Then by $(ii)$,
$[u_1]=[a_1[\overline{s_1}]b_1]$ where $[a_1s_1b_1]$ is a normal
$S$-diword, a contradiction.

This shows (iii).

$(iii)\Rightarrow (ii) $. \ Assume (iii).  Let $ 0\neq f\in Id(S)$.
Since the elements in $Irr(S)$ are linearly independent in $D(X|S)$,
by Lemma \ref{e4}, $[\bar{f}]=[a[\bar{s}]b]$, where $[asb]$ is a
normal $S$-diword. Thus, (ii) follows. \ \ $\square$

\ \

\noindent{\bf Remark:} In general, $(iii)\not\Rightarrow (i). $ For
example, it is noted that
$$
Irr(S)=\{x_j\dashv
x_{i_1}\dashv \dots \dashv x_{i_k} \ | \ j\in I, i_p\in I-I_0, \
1\leq p\leq k, \ i_1\leq \dots \leq i_k, \ k\geq 0\}
$$
is a linear basis of $D(X|S)$ in Theorem \ref{t2}. Let
$$
S_1=\{x_j\vdash x_i-x_i\dashv x_j+\{x_i,x_j\},\ x_t\dashv
x_{i_0},i,j,t\in I,i_0\in I_0\}.
$$
Then $Irr(S_1)=Irr(S)$ is a linear basis of $D(X|S)$. But in the
proof of Theorem \ref{t2}, we know that $S_1$ is not a
Gr\"{o}bner-Shirshov basis of $D(X|S)$.
\section{Applications}

In this section, we give  Gr\"{o}bner-Shirshov bases for the
universal enveloping dialgebra of a Leibniz algebra, the bar
extension of a dialgebra, the free product of two dialgebras, and
the Clifford dialgebra. By using our Theorem \ref{cd}, we obtain
some normal forms for  dialgebras mentioned the above.

\begin{definition}(\cite{Lo93})
A k-linear space $L$ equipped with bilinear multiplication $[,]$ is
called a Leibniz algebra if for any $a,b,c\in L$,
$$
[[a,b],c]=[[a,c],b]+[a,[b,c]]
$$
i.e., the Leibniz identity is valid in $L$.
\end{definition}

It is clear that if $(D,\dashv,\vdash)$ is a dialgebra then
$D^{(-)}=(D,[,])$ is a Leibniz algebra, where $[a,b]=a\dashv
b-b\vdash a$ for any $a,b\in D$.

\ \

If $f$ is a Leibniz polynomial in variables $X$, then by $f^{(-)}$
we mean a dialgebra polynomial in $X$ obtained from $f$ by
transformation $[a,b]\mapsto a\dashv b-b\vdash a$.

\begin{definition}
Let $L$ be a Leibniz algebra. A dialgebra $U(L)$ together with a
Leibniz homomorphism $\varepsilon: L\rightarrow U(L)$ is called the
universal enveloping dialgebra for $L$, if the following diagram
 commute:

\setlength {\unitlength}{1cm}
\begin{picture}(7, 3)
\put(4.2,2.3){\vector(1,0){1.7}} \put(4.1, 2.0){\vector(0,-1){1.3}}
\put(6,2.2){\vector(-1,-1){1.7}} \put(5.1,1.0){$\exists !f$}
 \put(3.9,0.2){$D$}
\put(6,2.2){$U(L)$} \put(3.9,2.2){$L$} \put(4.9, 2.4){$\varepsilon$}
\put(3.4,1.3){$\forall\delta$}
\end{picture}\\
where $D$ is a dialgebra, $\delta$ is a Leibniz homomorphism and
$f:U(L)\rightarrow D$ is a  dialgebra homomorphism such that
$f\varepsilon=\delta$ (i.e., $\varepsilon: L\rightarrow U(L)$ is a
universal arrow  in the sense of S. MacLane \cite{SM}, p55).

An equivalent definition is as follows: Let $L=Lei(X|S)$ is a
Leibniz algebra presented by generators $X$ and definition relations
$S$. Then $U(L)=D(X|S^{(-)})$ is the dialgebra with generators $X$
and definition relations $S^{(-)}=\{s^{(-)}| s\in S\}.$
\end{definition}
\begin{theorem}\label{t2}
Let $\mathcal{L}$ be a Leibniz algebra over a field $k$ with the
product $\{,\}$. Let $\mathcal{L}_0$ be the subspace of
$\mathcal{L}$ generated by the set $\{\{a,a\}, \{a,b\}+\{b,a\} \ | \
a,b\in \mathcal{L}\}$. Let $\{x_i|i\in I_0\}$ be a basis of
$\mathcal{L}_0$ and $X=\{x_i|i\in I\}$ a well ordered basis of
$\mathcal{L}$ such that $I_0\subseteq I$. Let $U(L)=D( X|x_i\dashv
x_j-x_j\vdash x_i- \{ x_i,x_j \})$ be the universal enveloping
dialgebra for $L$ and the ordering $<$ on $[X^*]$ as before. Then
\begin{enumerate}
\item[(i)] \ $D( X|x_i\dashv x_j-x_j\vdash x_i- \{
x_i,x_j \} )=D(X|S)$, where $S$ consists of the following
polynomials:
\begin{eqnarray*}
(a)&&f_{ji}=x_j\vdash x_i-x_i\dashv x_j+\{x_i,x_j\} \ \  \ \ \  \ \
\ \ \ \  \ \ \  \ \ \  \ \ \ \
 \  \  \ \ \ \  \  \  \ \ \ \  (i,j\in I)\\
(b)&&f_{ji\vdash t}=x_j\vdash x_i\vdash x_t-x_i\vdash x_j \vdash
x_t+\{x_i,x_j\}\vdash x_t \ \ \ \ \  \
 \  \  \ \ \ \ \ \  (i,j,t\in I, \ j>i)\\
(c)&&h_{i_0\vdash t}=x_{i_0}\vdash x_t \ \ \ \ \ \ \ \  \  \ \ \ \ \
\  \ \ \ \ \ \ \ \ \ \ \  \  \ \ \ \  \  \  \ \ \ \ \ \ \ \ \ \ \ \
\ \ \ \  \  \  \ \ \ \
(i_0\in I_0, \ t\in I)\\
(d)&&f_{t\dashv ji}=x_t\dashv x_j\dashv x_i-x_t\dashv x_i\dashv
x_j+x_t\dashv \{x_i,x_j \} \ \ \ \ \ \
 \  \  \ \ \ \ \ \  (i,j,t\in I, \ j>i)\\
(e)&&h_{t\dashv i_0}=x_t\dashv x_{i_0} \ \ \ \ \ \ \ \  \ \ \ \ \ \
\ \ \ \ \ \ \ \ \ \ \ \  \  \ \ \ \  \  \  \ \ \ \ \ \ \ \ \ \ \ \ \
\ \ \  \  \  \ \ \ \ (i_0\in I_0, \ t\in I)
\end{eqnarray*}
\item[(ii)] \ $S$ is a Gr\"{o}bner-Shirshov basis in $D(X)$.
\item[(iii)] \ The set
$$
\{x_j\dashv x_{i_1}\dashv \dots \dashv x_{i_k} \ | \ j\in I, i_p\in
I-I_0, \ 1\leq p\leq k, \ i_1\leq \dots \leq i_k, \ k\geq 0\}
$$
is a linear basis of the universal enveloping algebra
$U(\mathcal{L})$. In particular, $\mathcal{L}$ is a Leibniz
subalgebra of $U(\mathcal{L})$.
\end{enumerate}
\end{theorem}

\noindent{\bf Proof.}\ (i) By using the following

$$
f_{ji\vdash t}=f_{ji}\vdash x_t \
 \mbox{ and } \ f_{ji}\vdash x_t+f_{ij}\vdash x_t=(\{x_i,x_j\}+\{x_j,x_i\})\vdash x_t,
$$
we have (b) and (c) are in $Id(f_{ji})$. By symmetry, (d) and (e)
are in $Id(f_{ji})$. This shows (i).

(ii) We will prove that all compositions in $S$ are trivial modulo
$S$ (and $[w]$).  For convenience, we extend linearly the functions
$f_{ji}, \ f_{ji\vdash t}, \ f_{t\dashv ji}, \ h_{i_0\vdash t}$ and
$h_{t\dashv i_0}$ to $f_{j\{p,q\}} \ (f_{\{p,q\}i}), \ f_{ji\vdash
\{p,q\}}$ and $h_{\{p,q\}\dashv i_0}$, etc respectively. For
example, if $\{x_p,x_q\}=\sum\alpha_{pq}^sx_s$, then
\begin{eqnarray*}
f_{j\{p,q\}}&=&x_j\vdash\{x_p,x_q\}-\{x_p,x_q\}\dashv
x_j+\{\{x_p,x_q\},x_j\}=\sum\alpha_{pq}^sf_{js},\\
f_{ji\vdash\{p,q\}}&=&\sum\alpha_{pq}^s(x_j\vdash x_i\vdash
x_s-x_i\vdash x_j\vdash x_s+\{x_i,x_j\}\vdash
x_s)=f_{ji}\vdash\{x_p,x_q\},\\
h_{\{p,q\}\dashv i_0}&=&\sum\alpha_{pq}^sh_{s\dashv i_0}.
\end{eqnarray*}

By using the Leibniz identity,
\begin{equation}\label{e5}
\{\{a,b\},c\}=\{a,\{b,c\}\}+\{\{a,c\},b\},
\end{equation}
we have
\begin{eqnarray*}
\{a,\{b,b\}\}=0 \ \  and \ \ \{a,\{b,c\}+\{c,b\}\}=0
\end{eqnarray*}
for any $a,b,c\in \mathcal{L}$. It means that for any $i_0\in I_0, \
j\in I$,
\begin{equation}\label{e7}
\{x_j,x_{i_0}\}=0
\end{equation}
and by noting that
$\{x_{i_0},x_j\}=\{x_j,x_{i_0}\}+\{x_{i_0},x_j\}$, we have
\begin{equation}\label{e8}
\{x_{i_0},x_j\} \in \mathcal{L}_0.
\end{equation}
This implies that  $\mathcal{L}_0$ is an ideal of  $\mathcal{L}$.
Clearly, $\mathcal{L}/\mathcal{L}_0$ is a Lie algebra.

 The formulas (\ref{e5}),
(\ref{e7}) and (\ref{e8}) are useful in the sequel.

In $S$, all the compositions are as follows.

1) Compositions of left or right multiplication.

All possible compositions in $S$ of left multiplication are ones
related to (a), (b) and (c).

By noting that for any $s,i,j,t\in I$, we have
\begin{eqnarray*}
x_s\dashv f_{ji}&=&f_{s\dashv ji}  \ \ \ \ \ \  \  (j>i),\\
x_s\dashv f_{ji}&=&-f_{s\dashv
ij}+x_s\dashv(\{x_i,x_j\}+\{x_j,x_i\}) \ \ \ \ \  \ \  (j<i),\\
x_s\dashv f_{ii}&=&x_s\dashv\{x_i,x_i\}, \\
x_s\dashv f_{ji\vdash t}&=&f_{s\dashv ji}\dashv x_t \ \ (j>i) \ \ \ \ \ \mbox{and}\\
x_s\dashv h_{i_0\vdash t}&=&h_{s\dashv i_0}\dashv x_t,
\end{eqnarray*}
it is clear that all cases are trivial modulo $S$.

By symmetry, all compositions in $S$ of right multiplication are
trivial modulo $S$.

2) Compositions of inclusion and intersection.

We denote, for example, $(a\wedge b)$ the composition of the
polynomials of type $(a)$ and type $(b)$. It is noted that since (b)
and (c) are both left normed,  we have to prove that the
corresponding compositions of the cases of $(b\wedge b), \ (b\wedge
c), \ (c\wedge c)$ and $(c\wedge b)$ must be a linear combination of
left normed $S$-diwords in which the leading term of each $S$-diword
is less than $w$. Symmetrically, we consider the cases for the right
normed (d) and (e).

All possible compositions of inclusion and intersection are as
follows.

\begin{enumerate}
\item[($a\wedge c$)]\  $[w]=x_{i_0}\vdash
x_i \ (i_0\in I_0)$. We have, by (\ref{e7}),
\begin{eqnarray*}
&&(f_{i_0i},h_{i_0\vdash i})_{[w]} =-x_i\dashv
x_{i_0}+\{x_i,x_{i_0}\} =-h_{i\dashv {i_0}}\equiv0 \ \ mod(S,[w]).
\end{eqnarray*}

\item[($a\wedge d$)]\ $[w]=x_j\vdash
x_i\dashv x_q\dashv x_p \ \ (q>p)$. We have
\begin{eqnarray*}
&&(f_{ji},f_{i\dashv qp})_{[w]}\\
&=&-x_i\dashv x_j\dashv x_q\dashv x_p+\{x_i,x_j\}\dashv x_q\dashv
x_p +x_j\vdash x_i\dashv x_p\dashv x_p -x_j\vdash
x_i\dashv\{x_p,x_q\}\\
&=&-x_i\dashv f_{j\dashv qp}+f_{\{i,j\}\dashv qp}+f_{ji}\dashv
x_p\dashv x_q-f_{ji}\dashv\{x_p,x_q\}\\
&\equiv&0 \ \ mod(S,[w]).
\end{eqnarray*}

\item[$(a\wedge e)$]\  $[w]=x_j\vdash
x_{i}\dashv x_{i_0} \ (i_0\in I_0)$. We have
\begin{eqnarray*}
(f_{j{i}},h_{{i}\dashv {i_0}})_{[w]} =-x_{i}\dashv x_{j}\dashv
x_{i_0}+\{x_{i},x_j\}\dashv x_{i_0} =-x_i\dashv h_{j\dashv
i_0}+h_{\{i,j\}\dashv i_0}\equiv0 \ \ mod(S,[w]).
\end{eqnarray*}

\item[($b\wedge a$)]\ There are two cases to consider: $[w]=x_j\vdash
x_i\vdash x_t$ and $[w]=x_j\vdash x_i\vdash x_t\vdash x_p$.

For $[w]=x_j\vdash x_i\vdash x_t \ \ (j>i)$,\ by (\ref{e5}),\ we
have
\begin{eqnarray*}
(f_{ji\vdash t},f_{it})_{[w]}&=&-x_i\vdash x_j\vdash
x_t+\{x_i,x_j\}\vdash x_t+x_j\vdash
x_t\dashv x_i-x_j\vdash\{x_t,x_i\}\\
&=&-x_i\vdash f_{jt}+f_{\{i,j\}t}+f_{jt}\dashv
x_i-f_{j\{t,i\}}+f_{i\{t,j\}}-f_{it}\dashv x_j+f_{t\dashv ji}\\
&\equiv&0 \ \ mod(S,[w]).
\end{eqnarray*}

For $[w]=x_j\vdash x_i\vdash x_t\vdash x_p \ \ (j>i)$, we have
\begin{eqnarray*}
&&(f_{ji\vdash t},f_{tp})_{[w]}\\
&=&-x_i\vdash x_j\vdash x_t\vdash x_p+\{x_i,x_j\}\vdash x_t\vdash
x_p +x_j\vdash x_i\vdash x_p\dashv x_t -x_j\vdash
x_i\vdash\{x_p,x_t\}\\
&=&-x_i\vdash x_j\vdash f_{tp}+\{x_i,x_j\}\vdash f_{tp}+f_{ji\vdash
p}\dashv x_t-f_{ji\vdash \{p,t\}}\\
&\equiv&0 \ \ mod(S,[w]).
\end{eqnarray*}

\item[$(b\wedge b$)]\ There are two cases to consider:\ $[w]=x_j\vdash
x_i\vdash x_t\vdash x_s\vdash x_p$ and $[w]=x_j\vdash x_i\vdash
x_t\vdash x_p$.

For $[w]=x_j\vdash x_i\vdash x_t\vdash x_s\vdash x_p \ \ (j>i,t>s)$,
we have
\begin{eqnarray*}
&&(f_{ji\vdash t},f_{ts\vdash p})_{[w]}\\
&=&-x_i\vdash x_j\vdash x_t\vdash x_s\vdash x_p+\{x_i,x_j\}\vdash
x_t\vdash x_s\vdash x_p+x_j\vdash x_i\vdash x_s\vdash x_t\vdash
x_p\\
&&-x_j\vdash x_i\vdash\{x_s,x_t\}\vdash x_p\\
&=&-x_i\vdash x_j\vdash f_{ts\vdash p}+\{x_i,x_j\}\vdash f_{ts\vdash
p}+f_{ji\vdash s}\vdash x_t\vdash x_p-f_{ji\vdash\{s,t\}}\vdash
x_p\\
&\equiv&0 \ \ mod(S,[w])
\end{eqnarray*}
since it is a combination of left normed $S$-diwords in which the
leading term of each $S$-diword is less than $w$.

For $[w]=x_j\vdash x_i\vdash x_t\vdash x_p \ \ (j>i>t)$, suppose
that
$$
\{x_i,x_j\}=\sum_{m\in I_1}\alpha_{ij}^mx_m
+\alpha_{ij}^tx_t+\sum_{n\in I_2}\alpha_{ij}^nx_n \ (m<t<n).
$$
Denote
$$
B_{t\vdash\{i,j\}\vdash p}=x_t\vdash\{x_i,x_j\}\vdash
x_p-\{x_i,x_j\}\vdash x_t\vdash x_p-\{x_t,\{x_i,x_j\}\}\vdash x_p.
$$
Then
$$
B_{t\vdash\{i,j\}\vdash p}= \sum_{m\in I_1}\alpha_{ij}^m f_{tm\vdash
p}-\sum_{n\in I_2}\alpha_{ij}^n f_{nt\vdash p}-\sum_{q\in
I_0}\beta_qh_{q\vdash p}
$$
is a linear combination of left normed $S$-diwords of length 2 or 3,
where
$$
\sum_{q\in I_0}\beta_qx_{q}=\sum_{m\in
I_1}\alpha_{ij}^m(\{x_t,x_m\}+\{x_m,x_t\})+\alpha_{ij}^t\{x_t,x_t\}.
$$
Denote
$$\sum_{l\in
I_0}\gamma_lx_{l}=-(\{x_j,\{x_t,x_i\}\}+\{\{x_t,x_i\},x_j\})+(\{x_i,\{x_t,x_j\}\}+\{\{x_t,x_j\},x_i\}).$$
Now, by (\ref{e5}),\ we have
\begin{eqnarray*}
&&(f_{ji\vdash t},f_{it\vdash p})_{[w]}\\
&=&-x_i\vdash x_j\vdash x_t\vdash x_p+\{x_i,x_j\}\vdash x_t\vdash
x_p +x_j\vdash x_t\vdash x_i\vdash
x_p-x_j\vdash\{x_t,x_i\}\vdash x_p\\
&=&-x_i\vdash f_{jt\vdash p}-B_{t\vdash\{i,j\}\vdash p}+f_{jt\vdash
i}\vdash x_p-B_{j\vdash\{t,i\}\vdash p}
+\sum_{l\in
I_0}\gamma_lh_{{l}\vdash p}\\
&&+B_{i\vdash\{t,j\}\vdash p}-f_{it\vdash j}\vdash x_p+x_t\vdash
f_{ji\vdash p}\\
&\equiv&0 \ \ mod(S,[w])
\end{eqnarray*}
since it is a combination of left normed $S$-diwords in which the
leading term of each $S$-diword is less than $w$.

\item[$(b\wedge c)$]\ There are three cases to consider:\
$[w]=x_j\vdash x_{i_0}\vdash x_t \ (i_0\in I_0), \ [w]=x_{j_0}\vdash
x_i\vdash x_t \  (j_0\in I_0)$  and $[w]=x_j\vdash x_i\vdash
x_{t_0}\vdash x_n \ (t_0\in I_0)$.

Case 1.\ $[w]=x_j\vdash x_{i_0}\vdash x_t \ \ (j>i_0, \ i_0\in
I_0)$. By (\ref{e8}), we can assume that $\{x_{i_0},x_j\}=\sum_{l\in
I_0}\gamma_lx_{l}$. Then, we have
$$
(f_{ji_0\vdash t},h_{i_0\vdash t})_{[w]}=-x_{i_0}\vdash x_j\vdash
x_t+\{x_{i_0},x_j\}\vdash x_t=-h_{i_0\vdash j}\vdash x_t+\sum_{l\in
I_0}\gamma_lh_{l\vdash t}\equiv0 \ \ mod(S,[w]).
$$

Case 2.\  $[w]=x_{j_0}\vdash x_i\vdash x_t \ \ (j_0>i, j_0\in
I_0)$.\
 By
(\ref{e7}),\ we have
$$
(f_{j_0i\vdash t},h_{j_0\vdash i})_{[w]}=-x_i\vdash x_{j_0}\vdash
x_t+\{x_i,x_{j_0}\}\vdash x_t=-x_i\vdash h_{{j_0}\vdash t}\equiv0 \
\ mod(S,[w]).
$$

Case 3.\ $[w]=x_j\vdash x_i\vdash x_{t_0}\vdash x_n   \ \
(j>i,t_0\in I_0)$. We have
\begin{eqnarray*}
(f_{ji\vdash {t_0}},h_{{t_0}\vdash n})_{[w]} &=&-x_i\vdash x_j\vdash
x_{t_0}\vdash x_n+\{x_i,x_j\}\vdash
x_{t_0}\vdash x_n\\
&=&(- x_i\vdash x_j+\{x_i,x_j\})\vdash h_{{t_0}\vdash n}\\
&\equiv&0 \ \ mod(S,[w]) .
\end{eqnarray*}

\item[$(b\wedge d)$]\ $[w]=x_j\vdash x_i\vdash x_t\dashv x_q\dashv
x_p \ \ (j>i,q>p)$. We have
\begin{eqnarray*}
&&(f_{ji\vdash t},f_{t\dashv qp})_{[w]}\\
&=&-x_i\vdash x_j\vdash x_t\dashv x_q\dashv
x_p+\{x_i,x_j\}\vdash x_t\dashv x_q\dashv x_p\\
&& +x_j\vdash x_i\vdash x_t\dashv x_p\dashv x_q
-x_j\vdash x_i\vdash x_t\dashv\{x_p,x_q\}\\
&=&-x_i\vdash x_j\vdash f_{t\dashv qp}+\{x_i,x_j\}\vdash f_{t\dashv
qp}+f_{ji\vdash t}\dashv x_p\dashv x_q-f_{ji\vdash
t}\dashv\{x_p,x_q\}\\
&\equiv&0 \ \ mod(S,[w]).
\end{eqnarray*}

\item[($b\wedge e$)]\  $[w]=x_j\vdash x_i\vdash x_t\dashv x_{n_0} \ \ (j>i,n_0\in
I_0)$. We have
\begin{eqnarray*}
(f_{ji\vdash t},h_{t\dashv {n_0}})_{[w]} &=&-x_i\vdash x_j\vdash
x_t\dashv x_{n_0}+\{x_i,x_j\}\vdash
x_t\dashv x_{n_0}\\
&=&(-x_i\vdash x_j+\{x_i,x_j\})\vdash h_{t\dashv {n_0}}\\
&\equiv&0 \ \ mod(S,[w]).
\end{eqnarray*}

\item[$(c\wedge a)$]\ There are two cases to consider:\ $[w]=x_{n_0}\vdash x_t \ (n_0\in I_0)$
and $[w]=x_{n_0}\vdash x_t\vdash x_s \ (n_0\in I_0)$.

For  $[w]=x_{n_0}\vdash x_t \ (n_0\in I_0)$, we have
$$
(h_{{n_0}\vdash t},f_{{n_0}t})_{[w]} =x_t\dashv
x_{n_0}-\{x_t,x_{n_0}\} =h_{t\dashv{n_0}}\equiv0 \ \ mod(S,[w]).
$$

For $[w]=x_{n_0}\vdash x_t\vdash x_s \ (n_0\in I_0)$,\ we have
\begin{eqnarray*}
(h_{{n_0}\vdash t},f_{ts})_{[w]} =x_{n_0}\vdash x_s\dashv
x_t-x_{n_0}\vdash\{x_s,x_t\}=h_{{n_0}\vdash s}\dashv
x_t-h_{{n_0}\vdash \{s,t\}}\equiv0 \ \ mod(S,[w]).
\end{eqnarray*}

\item[$(c\wedge b)$]\ $[w]=x_{n_0}\vdash
x_t\vdash x_s\vdash x_p \ \ (t>s,n_0\in I_0)$.  We have
\begin{eqnarray*}
(h_{{n_0}\vdash t},f_{ts\vdash p})_{[w]}&=&x_{n_0}\vdash x_s\vdash
x_t\vdash x_p-x_{n_0}\vdash\{x_s,x_t\}\vdash
x_p\\
&=&h_{{n_0}\vdash s}\vdash x_t\vdash x_p-h_{{n_0}\vdash
\{s,t\}}\vdash x_p\\
&\equiv&0 \ \ mod(S,[w]).
\end{eqnarray*}

\item[$(c\wedge c)$]\  $[w]=x_{n_0}\vdash x_{t_0}\vdash x_r \ (n_0,t_0\in I_0)$.
We have
\begin{eqnarray*}
(h_{{n_0}\vdash {t_0}},h_{{t_0}\vdash r})_{[w]} =0.
\end{eqnarray*}

\item[$(c\wedge d)$]\ $[w]=x_{n_0}\vdash
x_t\dashv x_q\dashv x_p \ \ (q>p,n_0\in I_0)$. We have
\begin{eqnarray*}
(h_{{n_0}\vdash t},f_{t\dashv qp})_{[w]} &=&x_{n_0}\vdash x_t\dashv
x_p\dashv x_q-x_{n_0}\vdash
x_t\dashv\{x_p,x_q\}\\
&=&h_{{n_0}\vdash t}\dashv (x_p\dashv x_q-\{x_p,x_q\})\\
&\equiv&0 \ \ mod(S,[w]).
\end{eqnarray*}

\item[$(c\wedge e)$]\ $[w]=x_{n_0}\vdash
x_t\dashv x_{s_0} \ ({n_0},{s_0}\in{I_0})$.  We have
\begin{eqnarray*}
(h_{{n_0}\vdash t},h_{t\dashv{s_0}})_{[w]}=0.
\end{eqnarray*}
\end{enumerate}

Since $(d\wedge d)$,\ $(d\wedge e)$,\ $(e\wedge d)$,\ $(e\wedge e)$
are symmetric with $(b\wedge b)$,\ $(b\wedge c)$,\ $(c\wedge b)$,\
$(c\wedge c)$ respectively, they have the similar representations.
We omit the details.

So, we show that $S$ is a Gr\"{o}bner-Shirshov basis.

(iii) Clearly, the mentioned set is just the set $Irr(S)$. Now, the
results follow from Theorem \ref{cd}. \ \ $\square$

\ \

A Gr\"{o}bner-Shirshov basis $S$ is called reduced if $S$ is a monic
set and no monomial in any element of the basis  contains the
leading words of the other elements of the basis as subwords.

\ \

\noindent{\bf Remark:} Let the notation be in Theorem \ref{t2}. Let
$S^{red}$ consist of the following polynomials:
\begin{eqnarray*}
(a)&&f_{ji}=x_j\vdash x_i-x_i\dashv x_j+\{x_i,x_j\} \ \  \ \ \  \ \
\ \ \ \  \ \ \  \ \ \  \ \ \ \
 \  \  \ \ \ \  \  \  \ \ \ \  (i\in I, j\in I-I_0)\\
(b)&&f_{ji\vdash t}=x_j\vdash x_i\vdash x_t-x_i\vdash x_j \vdash
x_t+\{x_i,x_j\}\vdash x_t \ \ \ \ \  \
 \  \  \ \ \ \ (i,j\in I-I_0,\ j>i,t\in I)\\
(c)&&h_{i_0\vdash t}=x_{i_0}\vdash x_t \ \ \ \ \ \ \ \  \  \ \ \ \ \
\  \ \ \ \ \ \ \ \ \ \ \  \  \ \ \ \  \  \  \ \ \ \ \ \ \ \ \ \ \ \
\ \ \ \  \  \  \ \ \ \ \
(i_0\in I_0, \ t\in I)\\
(d)&&f_{t\dashv ji}=x_t\dashv x_j\dashv x_i-x_t\dashv x_i\dashv
x_j+x_t\dashv \{x_i,x_j \} \ \ \ \ \ \
 \  \  \ \ \ \ (i,j\in I-I_0,\ j>i,t\in I)\\
(e)&&h_{t\dashv i_0}=x_t\dashv x_{i_0} \ \ \ \ \ \ \ \  \ \ \ \ \ \
\ \ \ \ \ \ \ \ \ \ \ \  \  \ \ \ \  \  \  \ \ \ \ \ \ \ \ \ \ \ \ \
\ \ \  \  \  \ \ \ \ \ (i_0\in I_0, \ t\in I)
\end{eqnarray*}
Then $S^{red}$ is a reduced Gr\"{o}bner-Shirshov basis for $D(X|S)$.

We have the following corollary.

\begin{corollary}\ (\cite{A})\
Let the notation be as in Theorem \ref{t2}. Then as linear spaces,
$U(\mathcal{L})$ is isomorphic to $\mathcal{L}\otimes
U(\mathcal{L/L}_0)$, where $U(\mathcal{L/L}_0)$ is the universal
enveloping of the Lie algebra $\mathcal{L/L}_0$.
\end{corollary}
{\bf Proof.} Clearly, $\{x_j\ |\ j\in I-I_0\}$ is a $k$-basis of the
Lie algebra $\mathcal{L/L}_0$.  It is well known that the universal
enveloping $U(\mathcal{L/L}_0)$ of the Lie algebra $\mathcal{L/L}_0$
has a $k$-basis
$$
\{ x_{i_1}x_{i_2} \dots  x_{i_k} \ |\ i_1\leq \dots \leq i_k, \
i_p\in I-I_0, \ 1\leq p\leq k,  \ k\geq 0\}.
$$
By using (iii) in Theorem \ref{t2}, the result follows.  \ \
$\square$

\begin{definition}
Let $D$ be a dialgebra. An element $e\in D$ is called a bar unit of
$D$ if $e\vdash x=x\dashv e=x$ for any $x\in D$.
\end{definition}

\begin{theorem}\label{t9}
Each dialgebra has a bar unit extension.
\end{theorem}
{\bf Proof.} Let $(D,\vdash,\dashv)$ be an arbitrary dialgebra over
a field $k$ and $A$  the ideal of $D$ generated by the set
$\{a\dashv b-a\vdash b | \ a,b\in D\}$. Let $X_0=\{x_{i_0}|{i_0}\in
I_0\}$ be a $k$-basis of $A$ and $X=\{x_i|i\in I\}$ a well ordered
$k$-basis of $D$ such that $I_0\subseteq I$. Then $D$ has a
presentation by the multiplication table  $D=D(X|S)$, where
$S=\{x_i\vdash x_j-\{x_i\vdash x_j\},\ x_i\dashv x_j-\{x_i\dashv
x_j\},\ i,j\in I\}$, where $\{x_i\vdash x_j\}$ and $\{x_i\dashv
x_j\}$ are linear combinations of $x_t, t\in I$.

Let  $D_1=D(X\cup \{e\}|S_1)$, where $S_1=S\cup \{e\vdash y-y,\
y\dashv e-y,\ e\dashv x_{0},\ x_0\vdash e\ | \ y\in X\cup\{e\},
x_0\in X_0\}$. Then $D_1$ is a dialgebra with a bar unit $e$.

Denote
\begin{eqnarray*}
1.&&f_{i\vdash j}=x_i\vdash x_j-\{x_i\vdash x_j\},\\
2.&&f_{i\dashv j}=x_i\dashv x_j-\{x_i\dashv x_j\},\\
3.&&g_{e\vdash y}=e\vdash y-y,\\
4.&&g_{y\dashv e}=y\dashv e-y,\\
5.&&h_{x_{i_0}\vdash e}=x_{i_0}\vdash e,\\
6.&&h_{e\dashv x_{i_0}}=e\dashv x_{i_0},
\end{eqnarray*}
where $i,j\in I,\ i_0\in I_0,\ y\in X\cup\{e\}$.

We  show that $\{x_t\dashv x_{i_0}\}=0$ and $\{x_{i_0}\vdash
x_t\}=0$ for any  $t\in I,\ i_0\in I_0$.

Since $x_{i_0}\in A$, we have $x_{i_0}=\sum\alpha_i (c_if_id_i)$,
where $f_i=a_i\dashv b_i-a_i\vdash b_i$, $\alpha_i\in k$, $a_i,
b_i\in D$ and $c_i,d_i\in X^*$.

Since  $x_t\dashv (c_i(a_i\dashv b_i-a_i\vdash b_i)d_i)=0$, we have
$\{x_t\dashv \{c_i\{a_i\dashv b_i-a_i\vdash b_i \}d_i\}\}=0$ for
each $i$. Then $\{x_t\dashv x_{i_0}\}=0$.

By symmetry, we have $\{x_{i_0}\vdash x_t\}=0$.

To prove the theorem, by using our Theorem \ref{cd}, it suffices to
prove that with the ordering on $[(X\cup\{e\})^*]$ as before, where
$x<e,\ x\in X$, $S_1$ is a Gr\"{o}bner-Shirshov basis in $D(X\cup
\{e\})$. Now, we show that all compositions in $S_1$ are trivial.

All possible compositions of left and right multiplication are:
$z\dashv f_{i\vdash j}$, $z\dashv g_{e\vdash y}$, $z\dashv
h_{x_{i_0}\vdash e}$, $f_{i\dashv j}\vdash z$, $g_{y\dashv e}\vdash
z$, $h_{e\dashv x_{i_0}}\vdash z$, $z\in X\cup\{e\}$.

For $z\dashv f_{i\vdash j},\ z=x_t\in X$, since $(x_t\dashv
x_i)\dashv x_j=x_t\dashv(x_i\vdash x_j)$, we have $\{\{x_t\dashv
x_i\}\dashv x_j\}=\{x_t\dashv \{x_i\vdash x_j\}\}$ and
\begin{eqnarray*}
&&x_t\dashv f_{i\vdash j}\\
&=&x_t\dashv x_i\dashv x_j-x_t\dashv \{x_i\vdash x_j\}\\
&=&f_{t\dashv i}\dashv x_j+f_{\{t\dashv i\}\dashv j}-f_{t\dashv
\{i\vdash j\}}+\{\{x_t\dashv x_i\}\dashv x_j\}-\{x_t\dashv
\{x_i\vdash x_j\}\}\\
&=&f_{t\dashv i}\dashv x_j+f_{\{t\dashv i\}\dashv j}-f_{t\dashv
\{i\vdash j\}}\\
&\equiv&0\ mod(S_1).
\end{eqnarray*}

For $z\dashv f_{i\vdash j}, z=e$, let $\{x_i\dashv x_j\}-\{x_i\vdash
x_j\}=\sum\alpha_{i_0}x_{i_0}$. Then
\begin{eqnarray*} e\dashv
f_{i\vdash
j}&=&e\dashv x_i\dashv x_j-e\dashv \{x_i\vdash x_j\}\\
&=&e\dashv (x_i\dashv x_j-\{x_i\dashv x_j\})+e\dashv \{x_i\dashv
x_j\}-e\dashv \{x_i\vdash x_j\}\\
&=&e\dashv f_{i\dashv j}+\sum\alpha_{i_0}h_{e\dashv x_{i_0}}\\
&\equiv&0\ mod(S_1).
\end{eqnarray*}

For $z\dashv g_{e\vdash y}$, we have
\begin{eqnarray*}
z\dashv g_{e\vdash y}=z\dashv e\dashv y-z\dashv y=(z\dashv
e-z)\dashv y=g_{z\dashv e}\dashv y\equiv 0\ mod(S_1).
\end{eqnarray*}

For $z\dashv h_{x_{i_0}\vdash e}$, we have
\begin{eqnarray*}
z\dashv h_{x_{i_0}\vdash e}=z\dashv x_{i_0}\dashv e=z\dashv
g_{x_{i_0}\dashv e}+z\dashv x_{i_0}.
\end{eqnarray*}
It is clear that  $z\dashv x_{i_0}=h_{e\dashv x_{i_0}}$ if $z=e$ and
 $z\dashv x_{i_0}=x_t\dashv
x_{i_0}-\{x_t\dashv x_{i_0}\}=f_{t\dashv i_0}$ if $z=x_t\in X$,
since $\{x_t\dashv x_{i_0}\}=0$. This implies that $z\dashv
h_{x_{i_0}\vdash e}\equiv 0\ mod(S_1)$.

Thus we show that all compositions of left multiplication in $S_1$
are trivial modulo $S_1$. By symmetry, all compositions of right
multiplication in $S_1$ are trivial modulo $S_1$.

Now, all possible ambiguities $[w]$ of compositions of intersection
in $S_1$ are:

$1\wedge 1$,  $[x_ix_j\dot{x_t}]$;  $1\wedge 2$,
$[x_i\dot{x_j}x_t]$; $1\wedge 4$, $[x_i\dot x_je]$; $1\wedge 5$,
$[x_ix_{i_0}\dot e]$.

$2\wedge 2$,  $[\dot{x_i}x_jx_t]$; $2\wedge 4$, $[\dot x_ix_je]$.

$3\wedge 1$,  $[ex_i\dot x_j]$; $3\wedge 2$, $[e\dot x_ix_j]$;
$3\wedge 3$, $[ee\dot y]$; $3\wedge 4$, $[e\dot ye]$; $3\wedge 5$,
$[ex_{i_0}\dot e]$; $3\wedge 6$, $[e\dot ex_{i_0}]$.

$4\wedge 4$,  $[\dot yee]$; $4\wedge 6$,  $[\dot yex_{i_0}]$.

$5\wedge 3$,  $[x_{i_0} e\dot y]$; $5\wedge 4$, $[x_{i_0}\dot ee]$;
$5\wedge 6$,  $[x_{i_0}\dot ex_{j_0}]$.

$6\wedge 2$,  $[\dot ex_{i_0}x_j]$; $6\wedge 4$, $[\dot e
x_{i_0}e]$.

In the above, all $i,j,t\in I$, $i_0,j_0\in I_0$ and $y\in
X\cup\{e\}$.

There is no composition of inclusion in $S_1$.

We will show that all compositions of intersection in $S_1$ are
trivial. We check only the cases of $1\wedge 2$, $1\wedge 5$ and
$4\wedge 6$. Others can be similarly proved.

For $1\wedge 2$, $[w]=[x_i\dot{x_j}x_t]$, since $(x_i\vdash
x_j)\dashv x_t=x_i\vdash (x_j\dashv x_t)$, we have $\{\{x_i\vdash
x_j\}\dashv x_t\}=\{x_i\vdash \{x_j\dashv x_t\}\}$ and
\begin{eqnarray*}
(1\wedge 2)_{[w]}&=&-\{x_i\vdash x_j\}\dashv
x_t+x_i\vdash \{x_j\dashv x_t\}\\
&=&-f_{\{i\vdash j\}\dashv t}+f_{i\vdash \{j\dashv t\}}-
\{\{x_i\vdash x_j\}\dashv x_t\}+\{x_i\vdash \{x_j\dashv x_t\}\}\\
&=&-f_{\{i\vdash j\}\dashv t}+f_{i\vdash \{j\dashv t\}}\\
&\equiv&0\ mod(S_1,{[w]}).
\end{eqnarray*}

For $1\wedge 5$, $[w]=[x_ix_{i_0}\dot e]$, since  $x_i\vdash
x_{i_0}\in A$, we have $\{x_i\vdash
x_{i_0}\}=\sum\alpha_{j_0}x_{j_0}$ and
\begin{eqnarray*}
(1\wedge 5)_{[w]}=\{x_i\vdash x_{i_0}\}\vdash
e=\sum\alpha_{j_0}h_{x_{j_0}\vdash e}\equiv0\ mod(S_1,{[w]}).
\end{eqnarray*}

For $4\wedge 6$, $[w]=[\dot yex_{i_0}]$, we have
 $(4\wedge 6)_{[w]}=-h_{e\dashv x_{i_0}}$ if $y=e$ and  $
(4\wedge 6)_{[w]}=-f_{t\dashv i_0}$ if $y=x_t\in X$ since
$\{x_t\dashv x_{i_0}\}=0$. Then $(4\wedge 6)_{[w]}\equiv0\
mod(S_1,[w])$.

Then all the compositions in $S_1$ are trivial.

The proof is complete.  \ \ $\square$

\ \

\noindent{\bf Remark:} Let the notation be as in the proof of
Theorem \ref{t9}.  Let $D'=D(X\cup \{e_j\}_J|S')$ be a dialgebra,
where $S'=S\cup \{e_j\vdash y-y, y\dashv e_j-y,\ e_j\dashv x_{0},
x_0\vdash e_j\ | \ y\in X\cup\{e_j\}_J, x_0\in X_0,\ j\in J\}$. Let
$J$ be a well ordered set. Then with the ordering on $[(X\cup
\{e_j\}_J)^*]$ as before, where $x_i<e_j$ for all $i\in I,\ j\in J$,
by a similar proof of Theorem \ref{t9}, $S'$ is a
Gr\"{o}bner-Shirshov basis in $D(X\cup \{e_j\}_J)$. It follows from
Theorem \ref{cd} that $D$ can be embedded into the dialgebra $D'$
while $D'$ has bar units $\{e_j\}_J$.

\ \
\begin{definition}
Let $D_1,D_2$ be dialgebras over a field $k$. The dialgebra
$D_1*D_2$ with two dialgebra homomorphisms \  $\varepsilon_1:
D_1\rightarrow D_1*D_2$, $\varepsilon_2:D_2\rightarrow D_1*D_2$ is
called the free product of $D_1,D_2$, if the following diagram
commute:

\setlength {\unitlength}{1cm}
\begin{picture}(7, 3)
\put(4.2,2.5){\vector(1,0){1.7}}\put(6.7,2.3){\vector(0,-1){1.9}}
\put(9.2,2.5){\vector(-1,0){1.7}}\put(4.1, 2.3){\vector(1,-1){2}}
\put(9.2,2.3){\vector(-1,-1){2}}
 \put(6.5,0){$D$}\put(6.9,1.3){$\exists ! f$}
\put(6,2.4){$D_1*D_2$} \put(3.5,2.4){$D_1$} \put(9.5,2.4){$D_2$}
\put(4.9, 2.6){$\varepsilon_1$} \put(8.3,
2.6){$\varepsilon_2$}\put(4.5,0.9){$\forall\delta_1$}\put(8.7,0.9){$\forall\delta_2$}
\end{picture}\\
where $D$ is a dialgebra, $\delta_1,\delta_2$ are dialgebra
homomorphisms and $f:D_1*D_2\rightarrow D$ is a  dialgebra
homomorphism  such that
$f\varepsilon_1=\delta_1,f\varepsilon_2=\delta_2$ (i.e.,
$(\varepsilon_1,\varepsilon_2):(D_1,D_2)\rightarrow
(D_1*D_2,D_1*D_2)$ is a universal arrow in the sense of S. Maclane
\cite{SM}).

An equivalent definition is as follows: Let $D_i=D(X_i|S_i)$ be a
presentation  by generators and defining relations with $X_1\cap
X_2=\varnothing$, $i=1,2$. Then $D_1*D_2=D(X_1\cup X_2|S_1\cup
S_2)$.
\end{definition}

\ \

Let $(D_1,\vdash,\dashv)$, $(D_2,\vdash,\dashv)$ be two dialgebras
over a field $k$, $A_1$  the ideal of $D_1$ generated by the set
$\{a\dashv b-a\vdash b | \ a,b\in D_1\}$ and $A_2$  the ideal of
$D_2$ generated by the set $\{c\dashv d-c\vdash d | \ c,d\in D_2\}$.
Let $X_0=\{x_{i_0}|{i_0}\in I_0\}$ be a $k$-basis of $A_1$ and
$X=\{x_i|i\in I\}$ a well ordered $k$-basis of $D_1$ such that
$I_0\subseteq I$. Let $Y_0=\{y_{l_0}|{l_0}\in J_0\}$ be a $k$-basis
of $A_2$ and $Y=\{y_l|l\in J\}$ a well ordered $k$-basis of $D_2$
such that $J_0\subseteq J$. Then $D_1$ and $D_2$ have multiplication
tables:
\begin{eqnarray*}
&&D_1=D(X|S_1),\ \ \ S_1=\{x_i\vdash x_j-\{x_i\vdash x_j\},\
x_i\dashv x_j-\{x_i\dashv x_j\},\ i,j\in I\},\\
&&D_2=D(Y|S_2),\ \ \ S_2=\{y_l\vdash y_m-\{y_l\vdash y_m\},\
y_l\dashv y_m-\{y_l\dashv y_m\},\ l,m\in J\}.
\end{eqnarray*}

The free product $D_1*D_2$ of $D_1$ and $D_2$ is
$$
D_1*D_2=D(X\cup Y|S_1\cup S_2).
$$

We order $X\cup Y$ by $x_i<y_j$ for any $i\in I,j\in J$. Then we
have the following theorem.

\begin{theorem}\label{t10}
(i) $S$ is a Gr\"{o}bner-Shirshov basis of $D_1*D_2=D(X\cup
Y|S_1\cup S_2)$, where $S$ consists of the following relations:
\begin{eqnarray*}
1. &&f_{x_i\vdash x_j}=x_i\vdash x_j-\{x_i\vdash x_j\},\ \ \ \ i,j\in I ,\\
2. &&f_{x_i\dashv x_j}=x_i\dashv x_j-\{x_i\dashv x_j\},\ \ \ \ i,j\in I, \\
3. &&f_{y_l\vdash y_m}=y_l\vdash y_m-\{y_l\vdash y_m\},\ \ \ \ l,m\in J, \\
4. &&f_{y_l\dashv y_m}=y_l\dashv y_m-\{y_l\dashv y_m\},\ \ \ \ l,m\in J, \\
5.&&h_{x_{i_0}\vdash y_l}=x_{i_0}\vdash y_l,\ \ \ \ \ \ \ \ i_0\in I_0, l\in J, \\
6.&&h_{y_l\dashv x_{i_0}}=y_l\dashv x_{i_0},\ \ \ \ \ \ \ \ i_0\in I_0, l\in J,\\
7.&&h_{y_{l_0}\vdash x_i}=y_{l_0}\vdash x_i,\ \ \ \ \ \ \ \ i\in I, l_0\in J_0,\\
8.&&h_{x_i\dashv y_{l_0}}=x_i\dashv y_{l_0},\ \ \ \ \ \ \ \ i\in I,
l_0\in J_0.
\end{eqnarray*}

(ii) $Irr(S)$, which is a $k$-linear basis of $D_1*D_2$, consists of
all elements $z_{-m}\cdots z_{-1}\dot z_0z_1\cdots z_n$, where $
m,n\geq 0, z_0\in X\cup Y, z_i\in (X\setminus X_0) \cup (Y\setminus
Y_0),-m\leq i\leq n,i\neq 0, \mbox{neither }
\{z_j,z_{j+1}\}\subseteq X\ \mbox{nor } \{z_j,z_{j+1}\}\subseteq Y,
-m\leq j\leq n-1$.
\end{theorem}
{\bf Proof.} By the proof of Theorem \ref{t9}, we have $\{x_i\dashv
x_{i_0}\}=0, \ \{x_{i_0}\vdash x_i\}=0, \ \{y_l\dashv y_{l_0}\}=0$
and $\{y_{l_0}\vdash y_l\}=0$ for any  $i\in I,\ i_0\in I_0,\ l\in
J,\ l_0\in J_0$.

Firstly, we prove that $h_{y_l\dashv x_{i_0}}\in Id(S_1\cup S_2)$
for any $i_0\in I_0,\ l\in J$.

Since $y_l\dashv (c_i(\{a_i\dashv b_i\}-\{a_i\vdash b_i \})d_i)=
y_l\dashv (c_i((a_i\dashv b_i-\{a_i\dashv b_i \})-(a_i\vdash
b_i-\{a_i\vdash b_i \})d_i)\in Id(S_1\cup S_2)$, we have $y_l\dashv
\{c_i\{a_i\dashv b_i-a_i\vdash b_i \}d_i\}\in Id(S_1\cup S_2)$ for
all $i,l$. Then $h_{y_l\dashv x_{i_0}}\in Id(S_1\cup S_2)$.

Similarly, we have $h_{x_{i_0}\vdash y_l },\ h_{y_{l_0}\vdash x_i},\
h_{x_i\dashv y_{l_0}}\in Id(S_1\cup S_2)$ for any $i\in I,\ i_0\in
I_0,\ l\in J,\ l_0\in J_0$.

Secondly, we will show that all compositions in $S$ are trivial.

All possible compositions of left and right multiplication are:
$z\dashv f_{x_i\vdash x_j},\ z\dashv f_{y_l\vdash y_m},\ z\dashv
h_{x_{i_0}\vdash y_l},\ z\dashv h_{y_{l_0}\vdash x_i},\ f_{x_i\dashv
x_j}\vdash z,\ f_{y_l\dashv y_m}\vdash z,\ h_{y_l\dashv
x_{i_0}}\vdash z,\ h_{x_i\dashv y_{l_0}}\vdash z,$ where $z\in X\cup
Y$.

By a similar proof in Theorem \ref{t9}, all compositions of left and
right multiplication mentioned the above are trivial modulo $S$.

Now, all possible ambiguities $[w]$ of compositions of intersection
in $S$ are:
\begin{eqnarray*}
&&1\wedge 1,[x_ix_j\dot x_t];1\wedge 2,[x_i\dot x_jx_t];1\wedge 5,[x_ix_{i_0}\dot y_l];1\wedge 8,[x_i\dot x_jy_{l_0}].\\
&&2\wedge 2,[\dot x_ix_jx_t];2\wedge 8,[\dot x_ix_jy_{l_0}].\\
&&3\wedge 3,[y_ly_m\dot y_t];3\wedge 4,[y_l\dot y_my_t];3\wedge 6,[y_l\dot y_mx_{i_0}];3\wedge 7,[y_my_{l_0}\dot x_i].\\
&&4\wedge 4,[\dot y_ly_my_t];4\wedge 6,[\dot y_ly_mx_{i_0}].\\
&&5\wedge 3,[x_{i_0} y_l\dot y_t];5\wedge 4,[x_{i_0}\dot y_ly_t];
5\wedge 6,[x_{i_0}\dot y_lx_{j_0}];5\wedge 7,[x_{i_0}y_{l_0}\dot x_t].\\
&&6\wedge 2,[\dot y_lx_{i_0}x_t];6\wedge 8,[\dot y_mx_{i_0}y_{l_0}].\\
&&7\wedge 1,[y_{l_0}x_i\dot x_j];7\wedge 2,[y_{l_0}\dot x_ix_j];7\wedge 5,
[y_{l_0} x_{i_0}\dot y_m];7\wedge 8,[y_{l_0}\dot x_iy_{m_0}].\\
&&8\wedge 4,[\dot x_iy_{l_0}y_t];8\wedge 6,[\dot x_iy_{l_0}x_{i_0}].
\end{eqnarray*}
There is no composition of inclusion in $S$.

We will show that all compositions of intersection in $S$ are
trivial. We check only the cases of $1\wedge 5$ and $2\wedge 8$.
Others can be similarly proved.

For $1\wedge 5$, $[w]=[x_ix_{i_0}\dot y_l]$, let $\{x_i\vdash
x_{i_0}\}=\sum\alpha_{t_0}x_{t_0}$. Then
$$
(1\wedge 5)_{[w]}=-\{x_i\vdash x_{i_0}\}\vdash
y_l=-\sum\alpha_{t_0}h_{x_{t_0}\vdash y_l}\equiv 0\ \ mod(S,[w]).
$$

For $2\wedge 8$, $[w]=[\dot x_ix_jy_{l_0}]$, let $\{x_i\dashv
x_{j}\}=\sum\alpha_{t}x_{t}$. Then
$$
(2\wedge 8)_{[w]}=-\{x_i\dashv x_j\}\dashv
y_{l_0}=-\sum\alpha_{t}h_{x_t\dashv y_{l_0}}\equiv 0\ \ mod(S,[w]).
$$

Then all the compositions in $S$ are trivial. This show (i).

(ii) follows from our Theorem \ref{cd}. \ \ $\square$

\ \

\begin{definition}
Let $X=\{x_1,\dots ,x_n\}$ be a  set, $k$ a field of characteristic
$\neq 2$ and $(a_{ij})_{n\times n}$  a non-zero symmetric matrix
over $k$. Denote
$$
D(X\cup\{e\}\ |\ x_i\vdash x_j+x_j\dashv x_i-2a_{ij}e,\ e\vdash
y-y,\ y\dashv e-y,\ x_i,x_j\in X,\ y\in X\cup\{e\})
$$
by $C(n,f)$. Then $C(n,f)$ is called a Clifford dialgebra.
\end{definition}
We order $X\cup\{e\}$ by $x_1<\dots<x_n<e$.

\begin{theorem}\label{t4.5}
Let the notation be as the above. Then
\begin{enumerate}
\item[(i)] $S$ is a Gr\"{o}bner-Shirshov basis of Clifford dialgebra $C(n,f)$,
where $S$ consists of the following relations:
\begin{eqnarray*}
&&1.\ f_{x_ix_j}=x_i\vdash x_j+x_j\dashv x_i-2a_{ij}e,\\
&&2.\  g_{e\vdash y}=e\vdash y-y,\\
&&3.\  g_{y\dashv e}=y\dashv e-y,\\
&&4.\  f_{y\dashv x_ix_j}=y\dashv x_i\dashv x_j+y\dashv x_j\dashv x_i-2a_{ij}y,\ \ (i>j),\\
&&5.\  f_{y\dashv x_ix_i}=y\dashv x_i\dashv x_i-a_{ii}y,\\
&&6.\ f_{x_ix_j\vdash y}=x_i\vdash x_j\vdash y+x_j\vdash x_i\vdash
y-2a_{ij}y,\ \ (i>j),\\
&&7.\ f_{x_ix_i\vdash y}=x_i\vdash x_i\vdash y-a_{ii}y,\\
&&8.\ h_{x_ie}=x_i\vdash e-e\dashv x_i,
\end{eqnarray*}
where $x_i,x_j\in X, y\in X\cup\{e\}$.
\item[(ii)] A $k$-linear basis of $C(n,f)$ is a set of all elements of the form $\dot yx_{i_1}\cdots x_{ik}$,
where $y\in X\cup \{e\}$, $x_{ij}\in X$ and $i_1<i_2<\cdots <i_k\ \
(k\geq 0)$.
\end{enumerate}
\end{theorem}
{\bf Proof.} Let $S_1=\{f_{x_ix_j},g_{e\vdash y},g_{y\dashv e} \ |\
x_i,x_j\in X, y\in X\cup\{e\}\}.$

Firstly, we will show that $f_{y\dashv x_ix_j},\ f_{y\dashv
x_ix_i},\ f_{x_ix_j\vdash y},\  f_{x_ix_i\vdash y},\  h_{x_ie}\in
Id(S_1)$.

In fact, $f_{y\dashv x_ix_j}=y\dashv f_{x_ix_j}+2a_{ij}g_{y\dashv
e}$ implies $f_{y\dashv x_ix_j},\ f_{y\dashv x_ix_i}\in Id(S_1)$. By
symmetry, we have $f_{x_ix_j\vdash y},\  f_{x_ix_i\vdash y}\in
Id(S_1)$.

If there exists  $t$ such that $a_{it}\neq 0$, then
$$
2a_{it}h_{x_ie}=f_{x_ix_i\vdash x_t}-x_i\vdash f_{x_i\vdash
x_t}+f_{x_i\vdash x_t}\dashv x_i-f_{x_t\dashv x_ix_i}\in Id(S_1).
$$
Otherwise, $a_{it}= 0$ for any $t$. Since $(a_{ij})\neq 0$, there
exists  $j\neq i$ such that $a_{jt}\neq 0$ for some $t$. Then
\begin{eqnarray*}
&&2a_{jt}h_{x_ie}\\
&=&f_{x_ix_j\vdash x_t}-x_i\vdash f_{x_j\vdash x_t}-x_j\vdash
f_{x_i\vdash x_t}+f_{x_i\vdash x_t}\dashv x_j+f_{x_j\vdash
x_t}\dashv x_i-f_{x_t\dashv x_ix_j}\in Id(S_1).
\end{eqnarray*}

This shows that $h_{x_ie}\in Id(S_1).$

Secondly, we will show that all compositions in $S$ is trivial.

All possible compositions of left and right multiplication are:
$z\dashv f_{x_ix_j},\ z\dashv g_{e\vdash y},\ z\dashv
f_{x_ix_j\vdash y},\ z\dashv f_{x_ix_i\vdash y},\ z\dashv h_{x_ie},\
f_{x_ix_j}\vdash z,\ g_{y\dashv e}\vdash z,\ f_{y\dashv
x_ix_j}\vdash z,\ f_{y\dashv x_ix_i}\vdash z,\ h_{x_ie}\vdash z,$
where $z\in X\cup\{e\}$. We just check the cases of $f_{y\dashv
x_ix_j}\vdash z$ and $h_{x_ie}\vdash z$. Others can be similarly
proved.

For $f_{y\dashv x_ix_j}\vdash z$, we have
\begin{eqnarray*}
f_{y\dashv x_ix_j}\vdash z=y\vdash x_i\vdash x_j\vdash z+y\vdash
x_j\vdash x_i\vdash z-2a_{ij} y\vdash z=y\vdash f_{x_ix_j\vdash
z}\equiv0\ mod(S).
\end{eqnarray*}
For $h_{x_ie}\vdash z$,
\begin{eqnarray*}
h_{x_ie}\vdash z=x_i\vdash e\vdash z-e\vdash x_i\vdash z=x_i\vdash
g_{e\vdash z}-g_{e\vdash x_i}\vdash z\equiv0\ mod(S).
\end{eqnarray*}

Now, all possible ambiguities $[w]$ of compositions of intersection
in $S$ are:
\begin{eqnarray*}
&&1\wedge 3,[x_i\dot x_je];1\wedge 4,[x_i\dot x_jx_mx_n]\ (m>n);1\wedge 5,[x_i\dot x_jx_nx_n].\\
&&2\wedge 1,[ex_i\dot x_j];2\wedge 2,[ee\dot y];2\wedge 3,[e\dot
ye];2\wedge 4,[e\dot yx_ix_j]\ (i>j);\\
&&\ \ \ \ 2\wedge 5,[e\dot yx_ix_i];2\wedge 6, [ex_ix_j\dot y]\
(i>j);2\wedge 7, [ex_ix_i\dot y];
2\wedge 8,[ex_i\dot e].\\
&&3\wedge 3,[\dot yee];3\wedge 4,[\dot yex_ix_j]\ (i>j);3\wedge 5,[\dot yex_ix_i].\\
&&4\wedge 3,[\dot yx_ix_je]\ (i>j);4\wedge 4,[\dot
yx_ix_jx_mx_n]\ (i>j,m>n),[\dot yx_ix_jx_t]\ (i>j>t);\\
&&\ \ \ \ 4\wedge 5,[\dot yx_ix_jx_tx_t]\ (i>j),[\dot yx_ix_jx_j]\ (i>j).\\
&&5\wedge 3,[\dot yx_ix_ie];5\wedge 4,[\dot yx_ix_ix_mx_n]\
(m>n),[\dot yx_ix_ix_j]\ (i>j);\\
&&\ \ \ \ 5\wedge 5,[\dot yx_ix_ix_mx_m],[\dot yx_ix_ix_i].\\
&&6\wedge 1,[x_ix_jx_m\dot x_n]\ (i>j);6\wedge
2,[x_ix_je\dot y]\ (i>j); 6\wedge 3,[x_ix_j\dot ye]\ (i>j);  \\
&&\ \ \ \ 6\wedge 4,[x_ix_j\dot yx_mx_n]\ (i>j,m>n);6\wedge
5,[x_ix_j\dot yx_mx_m]\ (i>j);\\
&&\ \ \ \ 6\wedge 6,[x_ix_jx_mx_n\dot
y]\ (i>j,m>n),[x_ix_jx_t\dot y]\ (i>j>t); \\
&&\ \ \ \ 6\wedge 7,[x_ix_jx_mx_m\dot y]\ (i>j),[x_ix_jx_j\dot y]\
(i>j);
6\wedge 8,[x_ix_jx_t\dot e]\ (i>j).\\
&&7\wedge 1,[x_ix_ix_m\dot x_n];7\wedge 2,[x_ix_ie\dot y]; 7\wedge
3,[x_ix_i\dot ye]; 7\wedge 4,[x_ix_i\dot
yx_mx_n]\ (m>n); \\
&&\ \ \ \ 7\wedge 5,[x_ix_i\dot yx_mx_m];7\wedge 6,[x_ix_ix_mx_n\dot
y]\ (m>n),[x_ix_ix_t\dot y]\ (i>t);\\
&&\ \ \ \ 7\wedge 7,[x_ix_ix_mx_m\dot y],\ [x_ix_ix_i\dot y];7\wedge 8,[x_ix_ix_j\dot e].\\
&&8\wedge 3,[x_i\dot ee];8\wedge 4,[x_i\dot ex_mx_n]\ (m>n);8\wedge
5,[x_i\dot ex_mx_m].
\end{eqnarray*}

All possible ambiguities $[w]$ of compositions of  inclusion in $S$
are:
\begin{eqnarray*}
&&6\wedge 1,[x_ix_j\dot x_t]\ (i>j);\ 6\wedge 8,[x_ix_j\dot e]\
(i>j).\\
&&7\wedge 1,[x_ix_i\dot x_j];\ 7\wedge 8,[x_ix_i\dot e].
\end{eqnarray*}

 We just check the cases of intersection $1\wedge 4,
4\wedge 4, 6\wedge 4,6\wedge 8,8\wedge 4 $ and of inclusion $6\wedge
1,6\wedge 8$. Others can be similarly proved.

For $1\wedge 4$, $[w]=[x_i\dot x_jx_mx_n]\ (m>n)$, we have
\begin{eqnarray*}
&&(1\wedge 4)_{[w]}\\
&=& x_j\dashv x_i\dashv x_m\dashv x_n-2a_{ij}e
\dashv x_m\dashv x_n-x_i\vdash x_j\dashv x_n\dashv
x_m+2a_{mn}x_i\vdash x_j\\
&=&x_j\dashv f_{x_i\dashv x_mx_n}-2a_{ij}f_{e\dashv
x_mx_n}-f_{x_ix_j}\dashv x_n\dashv x_m+2a_{mn}f_{x_ix_j}\\
&\equiv& 0\ mod(S,[w]).
\end{eqnarray*}

For $4\wedge 4$, there are two cases to consider: $[w_1]=[\dot
yx_ix_jx_mx_n]\ (i>j,m>n)$ and $[w_2]=[\dot yx_ix_jx_t]\ (i>j>t)$.
We have
\begin{eqnarray*}
&&(4\wedge 4)_{[w_1]}\\
&=& y\dashv x_j\dashv x_i\dashv x_m\dashv
x_n-2a_{ij}y\dashv x_m\dashv x_n-y\dashv x_i\dashv x_j\dashv
x_n\dashv x_m+2a_{mn}y\dashv x_i\dashv x_j\\
&=&y\dashv x_j\dashv f_{x_i\dashv x_mx_n}-2a_{ij}f_{y\dashv
x_mx_n}-f_{y\dashv x_ix_j}\dashv x_n\dashv x_m+2a_{mn}f_{y\dashv
x_ix_j}\\
&\equiv& 0\ mod(S,[w_1])\ \ \ \ \ \mbox{and }\\
&& (4\wedge 4)_{[w_2]}\\
&=& y\dashv x_j\dashv x_i\dashv
x_t-2a_{ij}y\dashv x_t-y\dashv x_i\dashv x_t\dashv
x_j+2a_{jt}y\dashv x_i\\
&=& y\dashv f_{x_j\dashv x_ix_t}-f_{y\dashv x_jx_t}\dashv
x_i-f_{y\dashv x_ix_t}\dashv
x_j+y\dashv f_{x_t\dashv x_ix_j}\\
&\equiv& 0\ mod(S,[w_2]).
\end{eqnarray*}
For $6\wedge 4$, $[w]=[x_ix_j\dot yx_mx_n]\ (i>j,m>n)$, we have
\begin{eqnarray*}
&&(6\wedge 4)_{[w]}\\
&=&x_j\vdash x_i\vdash y\dashv x_m\dashv
x_n-2a_{ij}y\dashv x_m\dashv x_n-x_i\vdash x_j\vdash y\dashv
x_n\dashv x_m+2a_{mn}x_i\vdash x_j\vdash y\\
&=&x_j\vdash x_i\vdash f_{y\dashv x_mx_n}-2a_{ij}f_{y\dashv
x_mx_n}-f_{x_ix_j\vdash y}\dashv x_n\dashv
x_m+2a_{mn}f_{x_ix_j\vdash y}\\
 &\equiv& 0\ mod(S,[w]).
\end{eqnarray*}
For $6\wedge 8$, $[w]=[x_ix_jx_t\dot e]\ (i>j)$, we have
\begin{eqnarray*}
(6\wedge 8)_{[w]}
&=&x_j\vdash x_i\vdash x_t\vdash e-2a_{ij}x_t\vdash e+x_i\vdash x_j\vdash e\dashv x_t\\
&=&x_j\vdash x_i\vdash h_{x_te}-2a_{ij}h_{x_t
e}+f_{x_ix_j\vdash e}\dashv x_t\\
 &\equiv& 0\ mod(S,[w]).
\end{eqnarray*}
For $8\wedge 4$, $[w]=[x_i\dot ex_mx_n]\ (m>n)$, we have
\begin{eqnarray*}
(8\wedge 4)_{[w]}
&=&-e\dashv x_i\dashv x_m\dashv x_n-x_i\vdash e\dashv x_n\dashv x_m+2a_{mn}x_i\vdash e\\
&=&-e\dashv f_{x_i\dashv x_mx_n}-h_{x_i e}\dashv x_n\dashv
x_m+2a_{mn}h_{x_ie}\\
 &\equiv& 0\ mod(S,[w]).
\end{eqnarray*}

Now, we check the compositions of inclusion  $6\wedge 1$ and
$6\wedge 8$.

For $6\wedge 1$, $[w]=[x_ix_j\dot x_t]\ (i>j)$, we have
\begin{eqnarray*}
(6\wedge 1)_{[w]}
&=&x_j\vdash x_i\vdash x_t-2a_{ij}x_t-x_i\vdash x_t\dashv x_j+2a_{jt} x_i\vdash e\\
&=&x_j\vdash f_{x_i x_t}-f_{x_i x_t}\dashv x_j+2a_{jt} h_{x_i e}-f_{x_jx_t}\dashv x_i
+f_{x_t\dashv x_ix_j}+2a_{it}h_{x_je}\\
 &\equiv& 0\ mod(S,[w]).
\end{eqnarray*}
For $6\wedge 8$, $[w]=[x_ix_j\dot e]\ (i>j)$, we have
\begin{eqnarray*}
(6\wedge 8)_{[w]}
&=&x_j\vdash x_i\vdash e-2a_{ij}e+x_i\vdash e\dashv x_j\\
&=&x_j\vdash h_{x_i e}+h_{x_i e}\dashv x_j+h_{x_je}\dashv x_i+f_{e\dashv x_ix_j}\\
 &\equiv& 0\ mod(S,[w]).
\end{eqnarray*}

Then all the compositions in $S$ are trivial. We have proved (i).

For (ii),  since the mentioned set is just the set $Irr(S)$, by
Theorem \ref{cd} the result holds.

The proof is complete.  \ \ $\square$

\ \

\noindent{\bf Remark:} In the Theorem \ref{t4.5}, if the matrix
$(a_{ij})_{n\times n}=0$, then Clifford dialgebra $C(n,f)$ has a
Gr\"{o}bner-Shirshov basis $S'$ which consists of the relations
1--7.

 \ \

\noindent{\bf Acknowledgement}: The authors would like to thank P.S.
Kolesnikov who gives some valuable remarks for this paper.

\ \


\begin{thebibliography}{9}


\bibitem{A} M. Aymon and P.-P. Grivel, Un theoreme de Poincare-Birkhoff-Witt pour
les algebres de Leibniz, {\it Comm. Algebra}, 31(2003), N2, 527-544.


\bibitem{be78}G.M. Bergman, The diamond lemma for ring theory, {\it Adv. in
Math.}, 29, 178-218(1978).


\bibitem{bo72}L.A. Bokut, Unsolvability of the word problem, and
subalgebras of finitely presented Lie algebras, {\it Izv. Akad.
Nauk. SSSR Ser. Mat.}, 36, 1173-1219(1972).

\bibitem{bo76}L.A. Bokut, Imbeddings into simple associative
algebras, {\it Algebra i Logika}, 15, 117-142(1976).

\bibitem{bc1} L.A. Bokut and Yuqun Chen,  Gr\"{o}bner-Shirshov bases for Lie algebras: after
A.I. Shirshov,  {\it Southeast Asian Bull. Math.},  {31},
1057-1076(2007).

\bibitem{bo05} L.A. Bokut and K.P. Shum, Gr\"{o}bner and Gr\"{o}bner-Shirshov
bases in algebra: an elementary approach,  {\it Southeast Asian
Bull. Math.}, 29,  227-252(2005).

\bibitem{bu65}B. Buchberger, An algorithm for finding a basis for the
residue class ring of a zero-dimensional polynomial ideal [in
German], Ph.D. thesis, University of Innsbruck, Austria, (1965).

\bibitem{bu70}B. Buchberger, An algorithmical criteria for the
solvability of algebraic systems of equations[in German], {\it
Aequationes Math.}, 4, 374-383(1970).

\bibitem{K} P.S. Kolesnikov, Conformal representations of Leibniz algebras,
arXiv:math/0611501.

\bibitem{Lo93} J.-L. Loday, Une version non commutative des algebres de Lie:
les algebres de Leibniz, Ens. Math. 39, 269-293(1993).

\bibitem{Lo95} J.-L. Loday, Algebras with two associative operations (dialgebras),
C. R. Acad. Sci. Paris 321, 141-146(1995).


\bibitem{Lo99} J.-L. Loday, Dialgebras, in: Dialgebras and related operads,
Lecture Notes in Mathematics,
Vol. 1763. Berlin: Springer Verl., 2001, 7-66.

\bibitem{SM} S. MacLane, Categories for the Working Mathematician,
Springer, 1997.


\bibitem{sh62}A.I. Shirshov, Some algorithmic problem for Lie algebras,
 {\it Sibirsk. Mat. Z.}, 3(1962), 292-296(in Russian); English translation
in SIGSAM Bull., 33(2), 3-6(1999).




\end{thebibliography}
\end{document}